\documentclass[12pt]{amsart}

\usepackage{latexsym}
\usepackage[all]{xy}
\usepackage{amsmath, amsthm }
\usepackage{amssymb}
\usepackage{amsfonts}
\usepackage{color}
\usepackage{mathtools}
\usepackage[T1]{fontenc}
\usepackage{cancel}

\usepackage{hyperref}
\usepackage{xcolor}
\hypersetup{colorlinks=true,urlcolor=blue,linkcolor=blue,citecolor=blue}

\usepackage{cite}

\newcommand {\QCoh} {\mathbf{QCoh}}
\newcommand {\Coh} {\mathbf{Coh}}

\newcommand{\Gm}{\mathbb{G}_{\mathrm{m}}}

\newcommand {\Map} {\mathbf{Map}}
\newcommand {\Parf} {\mathbf{Perf}}
\newcommand {\rh} {\mathbb{R}\underline{Hom}}
\newcommand {\rch} {\mathbb{R}\underline{\mathcal{H}om}}
\newcommand {\OO} {\mathcal{O}}
\newcommand {\DR}{\mathsf{DR}}
\newcommand {\Fol}{\mathcal{F}ol}

\newcommand {\F} {\mathcal{F}}
\newcommand {\G} {\mathcal{G}}
\newcommand {\A} {\mathcal{A}}

\newcommand {\T} {\mathbb{T}}

\newcommand{\Z}{\mathbb{Z}}
\newcommand{\ZZ}{\mathbb{Z}}

\newcommand{\LL}{\mathbb{L}}
\newcommand{\VV}{\mathbb{V}}

\newcommand  {\dg}     {\mathbf{dg}}
\newcommand  {\fdg}     {\mathbf{dg}^{\mathrm{fil}}}
\newcommand  {\grdg}     {\mathbf{dg}^{\mathrm{gr}}}

\newcommand  {\edg}     {\epsilon-\mathbf{dg}}
\newcommand  {\medg}     {\epsilon\mathbf{dg}^{\mathrm{gr}}}

\newcommand  {\mecdga}     {\epsilon\mathbf{cdga}^{\mathrm{gr}}}

\newcommand{\s}{\infty}

\newcommand{\Alg}{\mathbf{Alg}}

\newcommand{\D}{\mathcal{D}}

\theoremstyle{plain}
\newtheorem{thm}{Theorem}[subsection]
\newtheorem{df}[thm]{Definition}
\newtheorem{prop}[thm]{Propositon}
\newtheorem{rmk}[thm]{Remark}
\newtheorem{cor}[thm]{Corollary}
\newtheorem{ex}[thm]{Example}
\newtheorem{lem}[thm]{Lemma}

\title[\resizebox{4.5in}{!}{Algebraic foliations and derived geometry II: the Grothendieck-Riemann-Roch theorem}]{Algebraic foliations and derived geometry: the Grothendieck-Riemann-Roch theorem}

\author{Bertrand To\"en and Gabriele Vezzosi}
\date{July 2020}

\begin{document}

\maketitle

\begin{abstract} 
This is the second of series of papers on the study of foliations in the setting of derived
algebraic geometry based on the central notion of \emph{derived foliations}. 
We introduce sheaf-like coefficients for derived foliations, called 
\emph{quasi-coherent crystals}, and construct a certain sheaf of dg-algebras of differential
operators along a given derived foliation, with the property that quasi-coherent crystals can be
interpreted as modules over this sheaf of differential operators. We use this interpretation in order
to introduce the notion of \emph{good filtrations} on quasi-coherent crystals, and define
the notion of characteristic cycle. Finally, we prove a Grothendieck-Riemann-Roch (GRR) formula expressing 
that formation of characteristic cycles is compatible with push-forwards
along proper and quasi-smooth morphisms. Several examples and applications are deduced from this, 
e.g. a GRR formula for $\D$-modules on possibly singular schemes.
\end{abstract}

\tableofcontents

\section*{Introduction}

This paper is a sequel to \cite{tvRH}, and the second in a series of works about derived techniques 
applied to the study of foliations. In \cite{tvRH}, we have introduced the notion of \emph{derived
foliations}, a far reaching generalization of the notion of foliation that is suited for the study of foliations with 
singularities\footnote{By the term \emph{singularities} of a foliation here we include two aspects: singularities of quotient type, such 
as foliations induced by a group scheme action (whose leaves are the orbits of the action), but also foliations by  singular subvarieties.}, and we have introduced the notion of locally free  crystals along a derived foliation, 
which, morally speaking, are perfect sheaves endowed with a flat connection along the leaves. We have proven,
under certain natural assumptions, that locally free crystals are in one-to-one correspondence with certain locally
constant sheaves of modules over the ring of flat functions, a correspondence that we have called the Riemann-Hilbert correspondence for derived foliations.\\

In the present work, we push further the study of categories of \emph{quasi-coherent} and \emph{coherent crystals} along derived foliations with the aim
of proving a Grothendieck-Riemann-Roch formula for them. We start by 
introducing \emph{rings of differential operators} along a derived foliations, generalizations of the usual rings of differential operators. In the derived context, these are not genuine sheaves of rings but rather sheaves of dg-algebras
with, in general, non-trivial cohomologies related to the singularities of the foliation. 
A first important result of the present work is that crystals can be interpreted
as sheaves of dg-modules over the ring of differential operators. \\

\medskip

\noindent \textbf{Theorem A}. Given a (possibly derived) scheme $X$ endowed with a derived foliation $\F$, 
there exists a sheaf of dg-algebras $\D_\F$ of differential operators along the leaves
of $\F$ such that there is an equivalence of $\s$-categories
$$\QCoh(\F) \simeq \D_\F-\dg_{X,qcoh}$$
between quasi-coherent crystals along $\F$ and sheaves of dg-modules over $\D_\F$ with 
quasi-coherent cohomologies. \\

\medskip

The above theorem must be understood in the spirit of the well known equivalence between crystals, defined
as sheaves on the crystalline site, and usual $\D$-modules (see for instance \cite[\S 5]{garo} for results in that
direction). We note also that when $\F=*_X$ is the final or tautological foliation (with a unique leaf given by $X$ istelf), 
the ring $\D_X$ should be called the ring of derived differential operators on $X$, and
the $\s$-category  $\D_X-\dg_{X,qcoh}$ already appears in \cite[Def. 4.2.9]{be}
under the name of category of \emph{derived $\D$-modules}. \\

We then study finiteness conditions on quasi-coherent crystals, by introducing the notion 
of \emph{coherent crystals}. They can be defined as compact objects in $\QCoh(\F)$, or as 
perfect sheaves of dg-modules over $\D_\F$. We study functorialities, such as 
pull-backs and push-forward along proper maps, and prove that proper and quasi-smooth push-forwards preserve
coherent crystals. This is a far reaching generalization of the well known fact that bounded coherent
complexes of $\D$-modules are stable by proper push-forward. 

One important consequence of Theorem A comes from the fact that 
$\D_\F$ is endowed with a canonical \emph{filtration}, induced by the degree of differential operators. It is this feature that allows
us to define the notion of a \emph{good filtration} on a coherent crystal, similar to the well-known  
notion for $\D$-modules. We could not convince ourselves that good filtrations always exist, 
but we prove they do exist for a large class of coherent crystals called \emph{finite cell} crystals, that  already contains
a lot of examples. A good filtration leads to the definition of \emph{characteristic cycle}
of a coherent crystal, by considering the associated graded object as a perfect complex
on the global cotangent stack $T^*\F \to X$ along the foliation. This characteristic cycle 
of a coherent crystal $E$,
is formally a $K$-theory class $Ch(E)\in K_0^{red}(T^*\F)$, in the reduced $K$-group of perfect
complexes on $T^*\F$, and is independent on the choice of the filtration. Here reduced K-groups are
defined as certain quotients of $K_0(T^*\F)$ by the classes of certain phantom objects. 
This is a technical point, related to the fact that working in the homotopical context 
implies existence of non-trivial good filtrations on the $0$-object, and we have to get rid
of them.\\

The second main result of this work is the following Grothendieck-Riemann-Roch formula
stating that taking characteristic cycles commute with proper quasi-smooth push-forwards.

\medskip

\noindent \textbf{Theorem B}.
Let $f : (X,\F) \longrightarrow (Y,\G)$ be a proper and quasi-smooth morphism between derived
schemes endowed with derived foliations. Then, for any coherent crystal $E$ along $\F$ that admits
a good filtration, we have
$$f_!(Ch(E))=Ch(f_!(E))$$
in $K^{red}_0(T^*\G)$.\\

\medskip

We deduce several applications of Theorem B, such as a Hirzebruch-Riemann-Roch formula for
\emph{foliated cohomology with coefficients} (see Corollary \ref{chrr}). We also present several other
examples of applications: an index formula for weakly Fredholm operators along a derived foliation (see Corollary \ref{cindex}), 
and a GRR formula for $\D$-modules on possibly singular schemes (\S \, \ref{s-grrdmod}). \\

The present work is organized in six sections. We start by relating filtered objects with 
graded mixed complexes, which is the core of the equivalence between crystals on the one side (defined
as graded mixed dg-modules) and $\D_\F$-modules on the other side. The second and third sections
are devoted to the general theory of crystals and sheaves of differential operators along derived foliations. 
Section 4 then defines basic functorialities, pull-backs and proper quasi-smooth push-forwards of crystals, 
and contains preservation properties of coherent objects. In Section 5, we introduce the global cotangent stack 
of a derived foliation, good filtrations and the notion of characteristic cycles. Finally, the last section 
contains the GRR formula, some ideas on how to extend the theorem to non-proper maps, as well as examples
and applications. \\

\noindent \textbf{Comments.} As a final remark in this introduction, we would like to emphasize that the results of this work are probably not optimal, in several different aspects. First of all, we are not completely satisfied with the existence and uniqueness statements concerning 
good filtrations for crystals on derived foliations. Another aspect is that we only consider quasi-coherent crystals, as opposed to the more general notion defined using Ind-coherent sheaves instead, as done in 
\cite{garo}. As a result, we only define push-forward for quasi-smooth morphisms, and we do not 
venture into defining the most general possible functorialities. This prevents us to 
obtain a full-fledged formalism of 6 operations for crystals on derived foliations, that, nonetheless,  
we are convinced should exist in a pretty general setting. We hope to be able to come back to this questions
in a future work. \\

\noindent \textbf{Acknowledgments.} This project has received funding from the European Research Council (ERC) under the European 
Union's Horizon 2020 research and innovation programme (grant agreement NEDAG ADG-741501).\\

\noindent\textbf{Notations and conventions.} We work over a base field $k$ of characteristic zero. 
All schemes, derived schemes and stacks are over $k$ and are assumed to be of finite presentation, i.e. locally of finite presentation, quasi-compact and quasi-separated. \\
%Unless otherwise specified, a functor will be understood as an $\s$-functor.

\section{Filtrations, mixed structures and filtered Tate realization}

\subsection{Filtered objects} In this short subsection we basically fix our notations on filtered objects and functors related to them. Everything here is gathered from existing literature (e.g. see \cite{tasos}).\\

\subsubsection{Generalities} Let $\mathcal{C}$ be a ($k$-linear) stable $\s$-category with sequential limits. We denote by $\mathcal{C}^{\mathrm{fil}}$ or by $\mathcal{C}^{\mathbb{Z}^{\leq}}$ the stable $\s$ category of $\s$-functors $\mathrm{Fun}(\mathbb{Z}^{\leq}, \mathcal{C})$ where 
$\ZZ^{\leq}$ is the ($\s$-nerve of the $1$-)category defined by the ordered set of integers (i.e. there is a unique map $i \to j$ iff $i\leq j$ in $\mathbb{Z}$). If $\mathbb{Z}$ denote the ($\s$-nerve of the $1$-) discrete category of integers (i.e. only identity maps are present), let us denote by $u:\mathbb{Z} \to \mathbb{Z}^{\leq}$ the obvious $\s$-functor. Restriction $u^*: \mathcal{C}^{\mathbb{Z}^{\leq}} \to \mathcal{C}^{\mathbb{Z}}$ along $u$ has a left adjoint $u_{!}:\mathcal{C}^{\mathbb{Z}} \to \mathcal{C}^{\mathbb{Z}^{\leq}}$, called the \emph{associated filtered object} $\s$-functor, given by left Kan extension along $u$, and $u_{!}$ has a left adjoint $\mathrm{Gr}=u^{!}:\mathcal{C}^{\mathbb{Z}^{\leq}} \to \mathcal{C}^{\mathbb{Z}}$ called the \emph{associated graded object} $\s$-functor. On objects, we have $$u_{!}((X_i)_{i\in \mathbb{Z}}) = \xymatrix{( \cdots \ar[r] & \coprod_{i\leq n-1}X_i \ar[r] & \coprod_{i\leq n}X_i \ar[r] & \cdots)}$$ and $$\mathrm{Gr}(Y)_n:=\mathrm{cofib}(Y_{n-1} \to Y_n)=: Y_n/Y_{n-1}.$$
We have natural $\s$-functors $(-)_{-\infty}:= \lim : \mathcal{C}^{\mathbb{Z}^{\leq}} \to \mathcal{C}$, and $(-)_{\infty}:= \mathrm{colim} : \mathcal{C}^{\mathbb{Z}^{\leq}} \to \mathcal{C}$. The functor $(-)_{\infty}$ is called the \emph{underlying object} $\s$-functor. 

If $\mathcal{C}$ is furthermore a \emph{presentable closed symmetric monoidal} stable category, then $\mathcal{C}^{\mathbb{Z}^{\leq}}$ inherits a symmetric monoidal structure (the so-called Day convolution) that can be described by $$(X\otimes Y)_n = \mathrm{colim}_{i+j \leq n} X_i \otimes_{\mathcal{C}} X_j.$$
For this tensor product, the underlying object $\s$-functor $(-)_{\infty}:= \mathrm{colim} : \mathcal{C}^{\mathbb{Z}^{\leq}} \to \mathcal{C}$ has a natural strong symmetric monoidal structure (i.e. the underlying object of the tensor product is naturally equivalent to the tensor product of the underlying objects).
The functor $\mathrm{Gr}:\mathcal{C}^{\mathbb{Z}^{\leq}} \to \mathcal{C}^{\mathbb{Z}}$ has a natural strong symmetric monoidal structure as well.\\

\subsubsection{Geometric picture}\label{sec-geompic} We now specialize, and describe geometrically, the previous notions in our case of interest 
where $\mathcal{C}$ is the $\s$-category $\dg :=\QCoh(\mathrm{Spec}\, k)$ of dg-$k$-modules. 

We let $\A:=[\mathbb{A}^1/\Gm]$ be the quotient stack of $\mathbb{A}^1$ by its natural
$\Gm$-action. 
The $\s$-category of (increasingly) filtered dg-modules is defined as 
$\fdg :=\QCoh(\A)$. By \cite[Thm. 1.1.]{tasos} we have an equivalence of symmetric monoidal $\s$-categories 
$\QCoh(\A)\simeq \mathrm{Fun}(\mathbb{Z}^{\leq}, \dg)$.
The $\s$-category 
of  graded dg-modules is defined as 
$\dg^{\mathrm{gr}} :=\QCoh(\mathrm{B}\mathbb{G}_{\mathrm{m}}) \simeq \mathrm{Fun}(\mathbb{Z},\dg)$.

Pull-back along the canonical map $1: \mathrm{Spec}\, k =\{1\} \to 
\mathbb{A}^1\to \mathbb{A}^1/\mathbb{G}_{\mathrm{m}}$ defines the \emph{underlying object} $\infty$-functor 
$(-)^u: \fdg \to \dg$ 
that can be viewed on objects as $(F^{\bullet}E):= \mathrm{colim}_i F^iE$.
Pull-back along the canonical map $0: \mathrm{B}\mathbb{G}_{m} \simeq \{0\}/\mathbb{G}_{m} \to 
\mathbb{A}^1/\mathbb{G}_{\mathrm{m}}$ 
defines the \emph{associated graded object} $\infty$-functor $\mathrm{Gr}: \fdg \to \dg^{\mathrm{gr}}$ that
can be viewed on objects as $\mathrm{Gr}(F^{\bullet}E):=\oplus_{i}F^i E/F^{i-1}E$.\\ Note that $\mathrm{Gr}$ 
is
left adjoint to the \emph{associated filtered object} $\s$-functor\footnote{This associated filtered object 
$\s$-functor 
can also be described as the left Kan extension functor along $\mathbb{Z} \to \mathbb{Z}^{\leq}$. It has also 
a right
adjoint $0^!$ equivalent to restriction along $\mathbb{Z} \to \mathbb{Z}^{\leq}$.} $0_*$ that, on objects, 
sends a 
graded dg-module $(V(i))_{i \in \mathbb{Z}}$ to the filtered dg-module $F^{\bullet}V$ where $F^p V:= 
\oplus_{i\leq p}V(i)$. 

The $\s$-category $\fdg$ comes equipped with a canonical symmetric monoidal structure, induced by 
tensor product of quasi-coherent complexes on the stack $\mathbb{A}^1/\mathbb{G}_{\mathrm{m}}$.\\

\subsubsection{Over a derived scheme $X$}\label{suX} Let now $X$ be a derived scheme over $k$. We will use the following notations:
\begin{itemize}
\item $\dg_X:= \mathsf{Sh}(X_{\mathrm{Zar}},\dg)$ the $\s$-category of (Zariski) sheaves of complexes of $k$-vector 
spaces. 
%\item $\medg_X :=\mathsf{Sh}(X_{\mathrm{Zar}},\medg)$ the $\s$-category of sheaves of mixed graded complexes 
%of $k$-vector spaces on $X_{\mathrm{Zar}}$.
\item $\grdg_X := \mathsf{Sh}(X_{\mathrm{Zar}},\grdg)$, the $\s$-category of (Zariski) sheaves of graded  complexes of 
$k$-vector spaces on $X_{\mathrm{Zar}}$.
\item $\fdg_X:= \mathsf{Sh}(X_{\mathrm{Zar}}, \fdg)$  the $\s$-category of (Zariski) sheaves of filtered complexes of 
$k$-vector spaces on $X_{\mathrm{Zar}}$.
\end{itemize}
The associated graded object functor $(-)^{\mathrm{gr}}: \fdg \to \grdg$ induces an associated graded object functor 
$(-)^{\mathrm{gr}}: \fdg_X \to \grdg_X$.   The forgetful functor $(-)^{\cancel{\epsilon}}: \medg \to \grdg$ induces a
forgetful functor $(-)^{\cancel{\epsilon}}: \medg_X \to \grdg_X$ forgetting the mixed structure. The 
underlying object functor $(-)^u: \fdg \to \dg$ induces an underlying object functor $(-)^u: \fdg_X \to 
\dg_X$.\\
Objects in $\dg_X$ will be referred to as \emph{complexes on $X$} (or 
equivalently \emph{over $X$}). Analogously, a \emph{filtered} (respectively, \emph{graded}) \emph{complex on} $X$ will be an object of $\fdg_X$ (respectively, of $\grdg_X$). Similar conventions will be adopted for ($E_1$-)algebras, and commutative algebras in $\dg_X$.

\subsection{Mixed structures}
We remind from \cite{cptvv} (see also the digest  \cite{pv}) the $\s$-category of graded mixed 
complexes (over $k$).
Its objects are $\Z$-graded objects $E=\oplus_{n\in \mathbb{Z}} E(n)$, inside the category of
cochain complexes  together
with extra differentials $\epsilon_n : E(n) \longrightarrow E(n+1)[-1]$, required to be morphisms of complexes.
These extra differentials  combine into
a morphism of graded complexes $\epsilon : E \longrightarrow E((1))[-1]$  
(where $E((1))$ is the graded complex obtained from $E$ by shifting the weight-grading
by $+1$), satisfying $\epsilon^2=0$. The datum of $\epsilon$ is called
a \emph{graded mixed structure} on the graded complex $E$. The complex $E(n)$ 
is itself called the \emph{weight n 
part} of $E$. 

Morphisms of
graded mixed complexes are defined in an obvious manner, and among them, the
quasi-isomorphisms are those morphisms inducing quasi-isomorphisms on all 
the weight-graded pieces \emph{individually}. 
By inverting quasi-isomorphisms, 
graded mixed complexes constitute an $\s$-category denoted by $\medg$. Alternatively, 
the $\s$-category $\medg$ can be defined as the $\s$-category of quasi-coherent
complexes $\QCoh(B\mathsf{H})$, over the classifying stack 
$B\mathsf{H}$ for the group stack $\mathsf{H}:=B\mathbb{G}_a \rtimes \mathbb{G}_m$
(see \cite[Rmk. 1.1.1]{cptvv} and \cite[Prop. 1.1]{pato}).

The $\s$-category $\medg$ comes equipped with a canonical symmetric monoidal 
structure. It is defined on objects by the usual 
tensor product of $\Z$-graded complexes (taken over the base field $k$), with the mixed structure
defined by the usual formula $\epsilon \otimes 1 + 1\otimes\epsilon$ (see \cite[\S 1.1]{cptvv}). 
When viewed as $\QCoh(B\mathsf{H})$, this is the usual
symmetric monoidal structure on quasi-coherent complexes on stacks. \\

\subsubsection{Over a derived scheme $X$}\label{s-medgX} Let now $X$ be a derived $k$-scheme. We will define 
$$\medg_X :=\mathsf{Sh}(X_{\mathrm{Zar}},\medg)$$ the symmetric monoidal $\s$-category of sheaves of mixed 
graded complexes of $k$-vector spaces on $X_{\mathrm{Zar}}$, the small Zariski site of $X$. 
Exactly as in \S \, \ref{suX}, objects in $\medg_X$ will be called \emph{graded mixed complexes on 
$X$}, and  we will allow ourselves to freely use expressions like  
graded mixed cdga's on $X$, filtered mixed complexes on $X$, filtered graded mixed dg-algebras on $X$, etc. The $\s$-category $\mathrm{CAlg}(\medg_X)$ of commutative algebras objects in $\medg_X$, i.e. of graded mixed cdga's over $X$, will be denoted  $\mecdga_X$.\\

\subsection{Filtered Tate realization} 
We start by introducing a symmetric lax monoidal $\s$-functor
$$|-|^t : \medg \longrightarrow \fdg,$$
called the \emph{filtered Tate realization} or simply the \emph{Tate realization}. This
construction already appears, in the non-filtered case, in \cite[\S 1.5]{cptvv}. 

For a graded mixed complex $E$, we define a filtered complex $|E|^t$ by the following formula
$$F^i|E|^t:=\prod_{p\geq -i}E(p)[-2p],$$
where we endow this infinite product with the usual total differential, sum of the
mixed structure and the cohomological differential. We have canonical 
inclusion maps
$F^i|E|^t \subset F^{i+1}|E|^t$, consisting of setting the first coordinate 
of $\prod_{p\geq -i-1}E(p)[-2p]$ to be zero. We thus have defined an
object in $\mathrm{Fun}(\ZZ^{\leq},\dg)$, and thus a filtered complex. By construction, the 
associated graded of $|E|^t$ is 
$\mathrm{Gr}(|E|^t) =\oplus_{q\in \mathbb{Z}} E(q)[-2q]$ where $E(q)[-2q]$ is of weight $-q$ 
(i.e. $\mathrm{Gr}^{q}F^{\bullet}:=F^q/F^{q-1}= E(-q)[2q]$ for any $q\in \mathbb{Z}$).
We warn the reader here concerning the change of signs between the weights of $E$ as a graded mixed
complex and the weights of $\mathrm{Gr}|E|^t$. 

Using the explicit formula above for $|-|^t$ we can identify the Tate realization as
a right adjoint in an adjunction  
$$k(*)\otimes_k - : \fdg \rightleftarrows \medg : |-|^t.$$
To see this we define an object 
$k(*)$ in $(\medg)^{\mathrm{fil}}$, i.e. a filtered graded mixed dg-module, by the formula
$$F^i k(*):=k(i)[2i],$$
where $k(i)$ is the trivial graded mixed complex pure of weight $i$. In order to define the 
transition morphism
$F^ik(*) \to F^{i+1}k(*)$, we must provide a morphism $k(i)[2i] \to k(i+1)[2i+2]$, or equivalently
a morphism $k \to k(1)[2]$, in $\medg$. Now we use that the mapping space 
$Map_{\medg}(k,k(1)[2])$ is discrete and canonically equivalent to the set $k$\footnote{In fact, $\mathbb{R}\underline{\mathrm{Hom}}_{\medg}(k,k(1)[2])$ (derived internal Hom in $\medg$) has all weights $n\neq 0$ complexes quasi-isomorphic to $0$, while its weight $0$ is quasi-isomorphic to $k[0]$. This can be computed, for example, using the cofibrant resolution $\tilde{k} \to k$ of $k$ in $\medg$ (in the projective model structure) of \cite[p. 503]{cptvv}. Simply notice the different conventions about mixed structures: while in \cite{cptvv} they have cohomological degree $1$, here they have degree $-1$. Use \cite[Rmk. 1.1.3]{cptvv} to get the $\tilde{k}$ needed here.}, and we take 
$1 \in k$ to define the required morphism $k(i)[2i] \to k(i+1)[2i+2]$. The filtered object
defined this way will be denoted by $k(*) \in \mathrm{Fun}(\ZZ^{\leq},\medg)$\footnote{Note that the pro-object
``lim''$_{i\leq 0} k(i)[2i]$ coincide with $k(-\s)$ appearing in \cite[\S 1.5]{cptvv}
(where a different convention for weights in graded mixed complexes was used)}.

The left adjoint of $|-|^t$ is then simply defined by tensoring a filtered 
complex with $k(*)$, using that the $\s$-category
of filtered objects in $\medg$ is naturally \emph{tensored} over $\fdg$. In a more explicit form, 
the left adjoint sends a filtered complex $E$ to the graded mixed complex
$$k(*)\otimes_k E := \int_{(p,q) \in \ZZ^{\leq}\times \ZZ^{\leq}}k(p)[2p]\otimes_k F^{-q}E,$$
defined as the coend of the $\s$-functor
$$\ZZ^{\leq}\times (\ZZ^{\leq})^{op} \longrightarrow \medg$$
sending $(p,q)$ to $k(p)[2p]\otimes_k F^{-q}E$. Note also that 
$|E|^t\simeq \rh(k(*),E)$, where we have implicitly used here the
equivalence $\mathrm{Fun}(\ZZ^{\leq},\dg)\simeq \mathrm{Fun}((\ZZ^{\leq})^{op},\dg)$ induced
by the isomorphism $p \mapsto -p$ between $\ZZ^{\leq}$ and $(\ZZ^{\leq})^{op}$
(which is the reason for the change of signs in weights).

The $\s$-functor $|-|^t$ clearly possesses a symmetric lax monoidal structure, 
coming from the fact that $k(*)$ has a natural structure of a commutative algebra 
inside filtered objects in $\medg$. This lax monoidal structure can also be seen
directly by the obvious explicit formulas on $F^i|E|^t=\prod_{p\geq -i}E(i)[-2i]$. 
Finally, also note that  
the underlying complex of the filtered complex $|E|^t$ is 
$$(|E|^t)^u= \mathrm{colim}_i \prod_{p\geq -i}E(p)[-2p] \in \dg.$$
To simplify notations, and when no confusion is possible, this underlying object will often be denoted simply by $|E|^t$.
Since $|-|^t$ 
is a symmetric lax monoidal $\s$-functor, it also induces $\s$-functors
on (associative) algebras and on their modules, a fact that we will use below. \\

We finish by the following proposition, stating that the Tate realization is
not very far from being an equivalence of $\s$-categories.

\begin{prop}\label{ptate}
The $\s$-functor
$$|-|^t : \medg \longrightarrow \fdg$$
is fully faithful, and its essential image consists of all 
filtered complexes $E$ that are complete, i.e. such that, for all $i\in \ZZ$, the natural morphism
$$F^iE \longrightarrow \lim_{j\leq i}(F^iE)/(F^jE)$$
is an equivalence.
\end{prop}

\noindent \textit{Proof.} For the fully faithfulness, we consider, for any $E \in \medg$, the counit map of the adjunction
$$k(*) \otimes |E|^t \longrightarrow E.$$
This is a morphism in the $\s$-category $\medg$. To check it is an equivalence we can forget the mixed structures
involved, as the forgetful $\s$-functor $\medg \to \grdg$ commutes with colimits (hence with coends).
Moreover, as a filtered object inside $\grdg$, 
$k(*)$ becomes the stupidly filtered object
$$\xymatrix{
\dots \ar[r] & k(i-1)[2i-2] \ar[r] & k(i)[2i] \ar[r] & k(i+1)[2i+2] \ar[r] & \dots
}$$
where all the maps $k(i)[2i] \to k(i+1)[2i+2]$ are zero. In other words, 
$k(*)$ can be written as a direct sum $\oplus 1(i)[2i]$ inside $\mathrm{Fun}(\ZZ^{\leq},\grdg)$, where
$1(i)$ is the filtered object in $\grdg$ with $F^i(1(i))=k(i)$ and $F^j(1(i))=0$ for all $j\neq i$.
Therefore, we have an equivalence of graded complexes
$$k(*)\otimes |E|^t \simeq \bigoplus_{i\in \ZZ}\mathrm{Gr}^{-i}(|E|^t)[2i],$$
where $\mathrm{Gr}^{-i}(|E|^t)[2i]$ sits in weight $i$.
Using this identification, we see that the counit morphism
$$k(*) \otimes |E|^t \longrightarrow E,$$
when considered as a morphism of graded complexes, is equivalent to the natural map
$$\oplus_i \mathrm{Gr}^{-i}(|E|^t)[2i] \longrightarrow \oplus_i E(i).$$
By the explicit formula for $|E|^t$ we see that this morphism is
indeed an equivalence of graded complexes. This shows that the counit of the adjunction $(k(*)\otimes-,|-|^t)$ is 
an equivalence, and thus that $|-|^t$ is fully faithful. \\
In order to characterize the essential image of $|-|^t$, we first 
notice, using the explicit formula for $|E|^t$, that the $\s$-functor $|-|^t$ does land inside the full sub-$\s$-category of \emph{complete} filtered complexes. To prove the statement it is thus enough
to show that the left adjoint $k(*)\otimes -$ is conservative when restricted to 
complete filtered complexes. For this we use the same argument as above: for an object 
$E \in \fdg$, the graded complex $k(*)\otimes E$ is equivalent to $\oplus_i \mathrm{Gr}^{-i}(E)[2i]$. 
But the functor $\mathrm{Gr} : \fdg \longrightarrow \grdg$ is obviously conservative on complete
objects, so this concludes the proof of the proposition.
\hfill $\Box$ \\

\subsubsection{Over a derived scheme $X$}\label{adjonX} If $X$ is a derived scheme, 
the above adjoint pair $(k(*)\otimes-,|-|^t)$ induces an analogous adjoint pair 
$$k(*)_X\otimes-: \fdg_X \rightleftarrows \medg_X: |-|_X^t, $$ 
and the obvious analog of Proposition \ref{ptate} holds for $|-|_X^t$. 
When no confusion is possible, the filtered Tate realization $|-|_X^t$ on $X$, will be again simply denoted as 
$|-|^t$.

\section{Quasi-coherent crystals on derived foliations}

In this section we introduce the notion of \emph{quasi-coherent crystals} along a derived
foliations (over a given derived scheme). We show that these can also be interpreted
as sheaves of modules over a certain sheaf of dg-algebras of differential operators (Theorem 
\ref{pdmod}), in the same manner as classical crystals over smooth varieties are equivalent to
left $\D$-modules. \\

\subsection{Quasi-coherent crystals}\label{sectionqcoh} 
Let us recall briefly from \cite{tvRH} the definition of a derived foliation on a derived scheme.

\begin{df}\label{d-dfol} Let $X$ be a derived scheme of finite presentation over $k$. The $\s$-category $\Fol(X)$ of \emph{derived foliations} on $X$ is the opposite of the full sub-$\s$-category of $\mecdga_X$ (\S \, \ref{s-medgX}) consisting of sheaves of graded mixed cdga's $\A$ such that  \begin{itemize}

\item $\A(0)\simeq \OO_X$.

\item The sheaf of $\OO_X$-dg-modules $\A(1)[-1]$ is perfect
and connective.

\item The natural morphism of sheaves of graded cdga's 
$$Sym_{\OO_X}(\A(1)) \longrightarrow \A^{\cancel{\epsilon}}$$
is a quasi-isomorphism.
\end{itemize}

For an object $\F \in \Fol(X)$, corresponding to a sheaf $\A$ of graded mixed cdga's on $X$, we will write $\LL_{\F}$ for $\A(1)[-1]$ (called the \emph{cotangent complex} of the foliation $\F$), and $\DR (\F)$ for $\A$ (called the \emph{de Rham algebra} of the foliation $\F$). \\ The \emph{initial} (respectively, \emph{final}) object in $\Fol(X)$, will be denoted by $0_X$ (respectively, $*_X$). Note that $\DR(0_X)= \OO_X$, while  $\DR(*_X)=\DR(X)$ is the derived de Rham algebra of $X$.

\end{df}

Let $\F \in \Fol(X)$ be a derived foliation on a 
derived scheme $X$, 
$\DR(\F) \in \mecdga_X=\mathbf{CAlg}(\medg_X)$ the corresponding sheaf of graded mixed cdga's on $X$, 
and $\DR(\F)-\medg_X :=\mathbf{Mod}_{\DR(\F)}(\medg_X)$ the 
$\s$-category of graded mixed $\DR(\F)$-dg-modules. 
By definition of derived foliation, $\DR(\F)^{\cancel{\epsilon}}\simeq Sym_{\OO_X}(\LL_{\F}[1])$ in $\mathbf{cdga}^{\mathrm{gr}}_{X}$,
so there is an induced morphism $\OO_X \to \DR(\F)^{\cancel{\epsilon}}$ in $\mathbf{cdga}^{\mathrm{gr}}_{X}$, where $\OO_X$ is concentrated in both weight and degree $0$. 

\begin{df}\label{defqcohcrys}
A \emph{quasi-coherent crystal over $\F$}
is a graded mixed $\DR(\F)$-dg-module $E$ satisfying the following two
conditions.

\begin{itemize}
\item The weight $0$ dg-module $E(0)$ is quasi-coherent over $\OO_X\simeq \DR(\F)(0)$.

\item The natural morphism
$$E(0)\otimes^{\LL}_{\OO_X}\DR(\F)^{\cancel{\epsilon}} \longrightarrow E$$
is a quasi-isomorphism of graded $\DR(\F)^{\cancel{\epsilon}}$-dg-modules over $X$.

\end{itemize}
The $\s$-category of quasi-coherent crystals over $\F$ is 
the full sub-$\s$-category $\QCoh(\F)$ of $\DR(\F)-\medg_X$
consisting of quasi-coherent crystals.
\end{df}

The $\s$-category $\QCoh(\F)$ is contravariantly functorial in the pair $(X,\F)$ in 
the following sense. Suppose that we have two pairs $(X,\F)$ and $(Y,\G)$ consisting
of derived schemes endowed with derived foliations. A morphism
$f : (X,\F) \longrightarrow (Y,\G)$ will consists of a pair
$(g,u)$, consisting of 
\begin{itemize}
\item a morphism $g : X \to Y$ of derived schemes,
\item a morphism $u : \F \to g^*(\G)$
of derived foliations over $X$ 
(i.e. a morphism $\DR(g^*\F):= g^*\DR(\G) \to \DR(\F)$ of graded mixed cdga's over $X$).
\end{itemize}

Associated to such a morphism there is a pull-back $\s$-functor
$$f^! : \QCoh(\G) \longrightarrow \QCoh(\F)$$
constructed as follows. By definition of pull-backs of derived foliations we have
an equivalence of graded mixed cdga's on $X$ 
$$\DR(g^*(\G))\simeq \DR(X) \otimes_{g^{-1}(\DR(Y))}g^{-1}(\DR(\G)).$$
The morphism $u$ thus corresponds to a morphism of graded mixed cdgas over $X$
under $\DR(X)$
$$u : \DR(X) \otimes_{g^{-1}(\DR(Y))}g^{-1}(\DR(\G)) \longrightarrow \DR(\F)$$
or equivalently, to a morphism of graded mixed cdga's under $g^{-1}(\DR(Y))$
$$u : g^{-1}(\DR(\G)) \longrightarrow \DR(\F).$$
The $\s$-functor $f^!$ on quasi-coherent crystals is thus simply defined by the following formula
(for $E$ a graded mixed $\DR(Y)$-module)
$$f^!(E):=g^{-1}(E) \otimes_{g^{-1}(\DR(\G))}\DR(\F).$$

Clearly, the rule $((X,\F) \mapsto \QCoh(\F), f \mapsto f^!)$ can be 
promoted to an $\s$-functor $\QCoh^{!}$ from the $\s$-category of pairs $(X,\F)$ to the
$\s$-category of $\s$-categories. Note also that 
$\QCoh(\F)$ comes equipped with a natural symmetric monoidal structure (induced by 
the tensor product of graded mixed $\DR(\F)$-modules), and that the pull-back $f^!$ has a natural symmetric
monoidal structure as well. 

Moreover, $f^!$ is compatible with the pull-back of quasi-coherent sheaves on derived schemes in the following sense.
We have a forgetful $\s$-functor
$$\QCoh(\F) \longrightarrow \QCoh(X)$$
which sends a graded mixed $\DR(\F)$-module $E$ to its weight zero part $E(0) \in \QCoh(X)$. For
a morphism $f=(g,u) : (X,\F) \longrightarrow (Y,\G)$ as above, the following square
is naturally commutative
$$\xymatrix{
\QCoh(\G) \ar[r]^-{f^!} \ar[d] & \QCoh(\F) \ar[d] \\
\QCoh(X) \ar[r]_-{g^*} & \QCoh(Y),
}$$
as this can be easily seen using the explicit formula
$f^!(E):=g^{-1}(E) \otimes_{g^{-1}(\DR(\G))}\DR(\F)$ and the condition stating that $E$ is 
a quasi-coherent crystal.

\subsection{Examples} We conclude this section by listing some examples $\QCoh(\F)$. \\

\noindent \textbf{Crystals over the trivial foliation.} When $\F=0_X$ is the initial foliation (often called, also, the trivial) foliation, defined by $\DR(\F):=\OO_X$ with trivial graded mixed structure, 
the $\s$-category $\QCoh(\F)$ clearly is equivalent to 
$\QCoh(X)$, the $\s$-category of quasi-coherent complexes over $X$. This
equivalence is realized by sending $E\in \QCoh(\F)$ to its weight $0$ part $E(0) \in \QCoh(X)$. It
can be promoted to an equivalence of symmetric monoidal $\s$-categories. \\

\noindent \textbf{Crystals over the tautological foliation.} Assume that $\F=*_X$ is the tautological 
foliation, that is the final object in $\Fol(X)$, defined by 
$\DR(\F):=\DR(X)$. If $X$ is a smooth variety, then there is a canonical equivalence of 
$\s$-categories
$$\QCoh(*_X)\simeq \D_X-\dg,$$
between quasi-coherent crystals along $*_X$ and 
quasi-coherent complexes of left $\D_X$-modules. This equivalence is constructed in 
\cite[\S 1.1]{pato}, but will be reviewed and generalized in the next section. This equivalence
is again compatible with the natural symmetric monoidal structures involved. \\

\noindent \textbf{Crystals over the Dolbeault foliation.} Let $\F=*_{Dol}$ be the Dolbeault foliation on 
a smooth variety $X$, defined by 
$\DR(\F):=Sym_{\OO_X}(\Omega_X^1[1]))$ where the graded cdga $Sym_{\OO_X}(\Omega_X^1[1]))$ is endowed with the \emph{zero} mixed structure. Then $\QCoh(*_{Dol})$ is naturally equivalent to 
the derived $\s$-category of complexes of quasi-coherent Higgs sheaves on $X$
(see \cite{sim}). \\

\noindent \textbf{Crystals over integrable foliations.} Let $f : X \longrightarrow Y$ be a morphism
of derived schemes and $\F:=f^*(0_Y)$ be the corresponding derived foliation (see \cite[1.3.3]{tvRH}).
Then $\QCoh(\F)$ is, by definition, the $\s$-category of \emph{relative
$\D$-modules on $X/Y$}. When $f$ is a smooth morphism between smooth varieties 
these relative $\D$-modules can be written as complexes of modules over $\D_{X/Y}$, the sheaf (of algebras) of
relative differential operators along the fibers of $f$. When $f$ is no more supposed to be smooth, 
we will see that $\D_{X/Y}$ only exists as a sheaf of \emph{dg-algebras} on $X$. 

\section{Rings of differential operators}

In this Section we associate to a derived foliation $\F$ on a derived scheme $X$, 
 a sheaf of filtered dg-algebras $\D_\F^{\mathrm{fil}}$ on $X$, called the sheaf of \emph{differential operators along} 
 $\F$, and prove that
 dg-modules over $\D_\F$ (filtration forgotten) corresponds to quasi-coherent modules along $\F$. \\

\subsection{Sheaf of differential operators} 

Let $X$ be a derived scheme and $\F$ be a derived foliation on $X$, with corresponding 
sheaf of graded mixed cdga's $\DR(\F)$, and  canonical augmentation 
$\DR(\F) \to \OO_X$ (making $\OO_X$ into a $\DR(\F)$-module). 

We consider $\rch_{\DR(\F)}(\OO_X,\OO_X)$, 
the sheaf of endomorphisms
of the graded mixed $\DR(\F)$-module $\OO_X$, i.e. the internal Hom object 
of endomorphisms of $\OO_X \in \DR (\F)-\medg_X$. This is a graded mixed dg-module over $\DR(\F)$, and a sheaf of
graded mixed $E_1$-algebra on $X$
\begin{equation}\label{Diff}
\rch_{\DR(\F)}(\OO_X,\OO_X) \in \mathbf{Alg}(\medg_X).
\end{equation} 
Its underlying graded dg-algebra, obtained by forgetting the mixed structure, is explicitly given
by
$$\rch_{\DR(\F)}(\OO_X,\OO_X)^{\cancel{\epsilon}} \simeq Sym_{\OO_X}(\mathbb{T}_\F[-2]),$$ 
where the tangent complex $\mathbb{T}_\F$ is here of weight $-1$. Note that, in 
general, the mixed structure induced on the right hand 
side is non-trivial and encodes invariants such as Atiyah classes, Lie brackets, etc.

\begin{df}
The \emph{filtered ring of differential operators of} $\F$ is 
defined to be 
$$\D_\F^{\mathrm{fil}}:=|\rch_{\DR(\F)}(\OO_X,\OO_X)|^t \in \mathbf{Alg}(\fdg_X).$$
The \emph{ring of differential operators along $\F$} is 
the underlying $E_1$-algebra over $X$ obtained by forgetting the filtration
and is denoted by
$$\D_\F:=(\D_\F^{\mathrm{fil}})^u \in \mathbf{Alg}(\dg_X).$$
\end{df}

By construction, $\D_\F^{\mathrm{fil}}$ is a sheaf of filtered ($k$-linear) dg-algebras on 
$X$. As usual, we set 
$$\D_\F^{\leq i}:=F^i\D_\F,$$
and call $\D_F^{\leq i}$ the \emph{sheaf of differential operators along $\F$ of order $\leq i$}.
Using the explicit description of the $\s$-functor $|-|^t$, it is straightforward 
to verify that the associated graded $\mathrm{Gr}(\D_\F^{\mathrm{fil}})$ is naturally equivalent to 
$Sym_{\OO_X}(\mathbb{T}_\F)$, where $\mathbb{T}_\F$ has pure weight $1$. \\

When $\F$ is a smooth derived foliation, and $X$ is a non-derived scheme (e.g. a smooth
variety), 
$\mathbb{T}_\F$ is a vector bundle on $X$, and thus
$\D_\F^{\mathrm{fil}}$ is automatically concentrated in degree $0$ i.e. it is a
genuine sheaf of filtered algebras over $X$. In general, $\D_\F$ is a dg-algebra on $X$,
that might have non-trivial cohomologies in an infinite number of degrees, and is bounded on the right. It is
moreover concentrated in non-negative degrees when $X$ is a non-derived scheme. 
When $\F$ and $X$ are both quasi-smooth (i.e. $\LL_\F$ and $\LL_X$ are perfect of tor-amplitude $[-1,0]$),
the dg-algebra $\D_F$ is moreover cohomologically bounded. This can be checked for instance using the exact triangles 
$$\xymatrix{\D_\F^{\leq i} \ar[r] & \D_\F^{\leq i+1} \ar[r] & Sym^{i+1}_{\OO_X}(\T_\F)
}$$
and induction on $i$. More is true: if $X$ is a non-derived scheme and 
$\T_\F$ can be represented by a two term complex of vector bundles $V \to W$, then $\D_\F$ is
cohomologically concentrated in degrees $[0,rk(W)]$. \\

Here are some basic examples of rings of differential operators.

\begin{ex}\label{Dex}
\begin{enumerate}
\emph{\item When $\F$ is the zero foliation (i.e. $\mathbb{L}_\F=0$)
then $\D_\F=\OO_X$ with the trivial filtration.
\item When $\F$ is the tautological foliation on a smooth variety $X$
(i.e. $\DR(\F)=\DR(X)$), then $\D_\F=\D_X$ is the usual ring
of differential operators with its usual filtration by the order of
operators.
\item When $\F$ is the Dolbeault foliation on a smooth variety, that 
is $\DR(\F)=Sym_{\OO_X}(\Omega^1_X[1])$ with trivial mixed structure, 
then $\D_\F=Sym_{\OO_X}(\mathbb{T}_X)$ with the split filtration. More generally,
when $\DR(\F)=Sym_{\OO_X}(\LL_\F[1])$ with trivial mixed structure (abelian 
derived foliation), then $\D_\F\simeq Sym_{\OO_X}(\T_\F)$. 
\item If the foliation $\F$ is smooth, i.e a Lie algebroid (see \cite{tvRH}), induced
by a smooth groupoid $G$ acting on $X$, then 
$\D_\F$ is the ring of distributions on $G$, i.e. the $\OO_X$-linear 
dual of formal functions of $G$, endowed with the convolution product.
This coincides with the universal enveloping algebra of the Lie algebroid.
\item When $\F$ is globally integrable by a flat 
and generically smooth
morphism of smooth varieties $f : X \longrightarrow Y$, 
$\D_\F$ is called the dg-algebra of \emph{relative differential operators}. 
The reason for this name comes from the fact that $H^0(\D_\F)$ is indeed a subring of 
$\D_X$ consisting of differential operators stabilizing the
fibers of $f$. }
\end{enumerate}
\end{ex}

\begin{rmk}\label{weyl2} \emph{It is interesting to note the following basic example. Let $X=\mathrm{Spec}( k[u])$ with $u$ in degree $-2$, and $\F=*_X$ be the tautological (i.e. final) derived foliation on $X$. Then
$\D_\F$ is here the Weyl dg-algebra over one generator $u$ in degree $-2$. In other
words, it is the dg-algebra freely generated by two cocycle $u$ and $\frac{\partial}{\partial u}$, 
respectively in degrees $-2$ and $2$, with the usual commutation relation
$[u,\frac{\partial}{\partial u}]=-1.$ Note that $\D_\F$ is \emph{not} Morita equivalent to $k$, as opposed to the case
when $\deg u$ is odd (see \cite[Proof of Cor. 4.3.13]{be}). The reader will find more about $\D_{*_X}$ in \cite{be}, even for more general $X$ than the ones we consider in this paper.}
\end{rmk}

The rule $\F \mapsto \D_{\F}$ defines an $\s$-functor $\Fol(X) \to \mathbf{Alg}(\fdg_X)$: if $\F \to \F'$ 
is a morphism in $\Fol (X)$, i.e. a morphism $\varphi: \DR(\F') \to \DR(\F)$ of sheaves of graded mixed 
cdga's over $X$, $\varphi$ induces a map $\rch_{\DR(\F)}(\OO_X,\OO_X) \to \rch_{\DR(\F')}(\OO_X,\OO_X)$,
thus an induced map $\D_\F^{\mathrm{fil}} \to \D_{\F'}^{\mathrm{fil}}$ in $\mathbf{Alg}(\fdg_X)$. In particular,
the maps $0_X \to \F \to *_X$, from the initial and to the final foliations, provides
maps of filtered dg-algebras over $X$
$$\xymatrix{
\OO_X \ar[r] & \D_\F^{\mathrm{fil}} \ar[r] & \D_X^{\mathrm{fil}},}$$
where $\OO_X$ has the trivial filtration, and $\D_X^{\mathrm{fil}}$ is by definition
the ring of differential operators on $X$ with its natural filtration by order of
operators. 
%(\textcolor{blue}{Ref to Dario ?}). \\

\subsection{Quasi coherent $\F$-crystals and $\D_\F$-modules}

As we have seen in Section \ref{adjonX}, the Tate realization provides a fully faithful 
symmetric lax monoidal $\s$-functor
$$|-|^t : \medg_X \longrightarrow \fdg_X.$$
Note that this full embedding is not symmetric monoidal, as the 
tensor product in $\fdg_X$ does not preserve complete filtered objects in general. However, 
$|-|^t$ sends the tensor product of $\medg_X$ to the completed tensor product in $\fdg_X$, where
the completion is taken with respect of the filtration.

We will now construct an $\s$-functor
$$\psi_X : \QCoh(\F) \longrightarrow \D_\F-\dg_X,$$
from quasi-coherent crystals along $\F$ to (left) $\D_\F$-dg-modules on $X$. We start, as in formula (\ref{Diff}) above, by letting 
$$S(\F):=\rch_{\DR(\F)}(\OO_X,\OO_X) \in \Alg(\medg_X),$$
the internal Hom object 
of endomorphisms of $\OO_X \in \mathrm{Mod}_{\DR (\F)}(\medg_X)$.
%We have already seen that 
%$S(\F)$ is equivalent, as a graded dg-algebra to $Sym_{\OO_X}(\T_\F[-2])$ i.e. 
%$$S(\F)^{\cancel{\epsilon}}\simeq Sym_{\OO_X}(\T_\F[-2])$$
%where $\mathbb{T}_\F$ is here in weight $-1$.
The object $\OO_X$ can thus be considered as a graded mixed bi-module with 
its right action by $\DR(\F)$ and left action by $S(\F)$. Therefore it can be used in order
to produce the following $\s$-functor (between categories of left modules)
$$\OO_X \otimes_{\DR(\F)} - : \DR(\F)-\medg_X \longrightarrow S(\F)-\medg_X$$
Note that this $\s$-functor obviously sends  
$\DR(\F)$-dg-modules which are graded free on weight zero, i.e. quasi-coherent crystals, to 
graded mixed $S(\F)$-dg-modules whose underlying graded module are pure of weight $0$. We thus get 
an induced $\s$-functor
$$\psi'_X: \QCoh(\F) \longrightarrow S(\F)-\medg_{qcoh,w=0},$$
where $S(\F)-\medg_{qcoh,w=0} \subset S(\F)-\medg_X$ is the full sub-$\s$-category of objects
whose underlying graded $S(\F)$-module are pure of weight $0$ and quasi-coherent over 
$\OO_X$.

We now compose the previous $\s$-functor with the Tate realization $|-|^t$ in order to get 
an $\s$-functor
$$\psi_X : \xymatrix{\QCoh(\F) \ar[r]^-{\psi'_X} & S(\F)-\medg_{qcoh,w=0} \ar[r]^-{|-|^t} & \D_\F-\dg.}$$ 

\begin{thm}\label{pdmod}
The $\s$-functor $\psi_X: \QCoh(\F) \to \D_\F-\dg$ defined above is fully faithful. Its essential image, $\D_\F-\dg^{qcoh}_{X}$, consists of 
all $\D_\F$-modules over $X$ which are quasi-coherent as $\OO_X$-modules
$$\psi_X : \QCoh(\F) \simeq \D_\F-\dg^{qcoh}_{X}.$$
\end{thm}

\noindent \textit{Proof.} We will first prove that the $\s$-functor
$$\psi'_X: \QCoh(\F) \longrightarrow S(\F)-\medg_{qcoh,w=0}$$
is an equivalence of $\s$-categories, and then identify the essential image of $\psi_X$. \\

\noindent \textbf{$\psi'_X$ is fully faithful.} To establish the fully faithfulness of $\psi'_X$, we start by noticing that the
induced morphism of graded mixed dg-algebras over $X$
$$\DR(\F)=\rch_{\DR(\F)}(\DR(\F),\DR(\F)) \longrightarrow \rch_{S(\F)}(\OO_X,\OO_X)$$
is an equivalence, as this can be checked directly at the graded algebras level by forgetting the mixed
structures. This easily implies that for $E,E' \in \QCoh(\F)$, the induced morphism of graded mixed
complexes over $X$
$$\rch_{\DR(\F)}(E,E') \longrightarrow \rch_{S(\F)}(\psi_X(E),\psi_X(E'))$$
is also an equivalence of graded mixed complexes over $X$. By passing to realizations $|-|$ on both
sides we get that the induced morphism of mapping spaces
$$\Map(E,E')\simeq |\rch_{\DR(\F)}(E,E')|  \longrightarrow \Map(\psi(E),\psi(E')) \simeq  
|\rch_{S(\F)}(\psi_X(E),\psi_X(E'))|$$
is an equivalence.\\

\noindent \textbf{$\psi'_X$ is essentially surjective.} To prove the essential surjectivity of $\psi'_X$ we use the $\s$-functor
$$\rch_{S(\F)}(\OO_X,-) : S(\F)-\medg \longrightarrow \DR(\F)-\medg$$
right adjoint to $\OO_X \otimes_{\DR(\F)}-$. Assume that $E \in S(\F)-\medg_{qcoh,w=0}$. 
Then, as a graded $S(\F)$-module, $E$ is induced from a quasi-coherent module
$E(0)$ over $\OO_X$ via the augmentation $S(\F) \to \OO_X$. As a consequence, 
the graded module underlying $\rch_{S(\F)}(\OO_X,E)$ is of the form 
$\DR(\F)\otimes_{\OO_X}E(0)$. In particular, $\rch_{S(\F)}(\OO_X,-)$
produces a right adjoint of the $\s$-functor $\psi_X$ restricted
to the sub-$\s$-categories under consideration
$$\OO_X \otimes_{\DR(\F)}- : \QCoh(\F) \rightleftarrows S(\F)-\medg_{qcoh,w=0} : \rch_{S(\F)}(\OO_X,-).$$
We already know that the left adjoint of the above adjunction is fully faithful. Moreover, 
the right adjoint is easily seen to be conservative, by using the above comments concerning the 
underlying graded objects. We thus conclude that $\psi'_X: \OO_X \otimes_{\DR(\F)}-$ indeed produces an 
equivalence of $\s$-categories
$$\QCoh(\F) \simeq S(\F)-\medg_{qcoh,w=0}.$$

\noindent \textbf{End of the proof.} To finish the proof of the theorem, we now consider the filtered Tate realization
$$|-|^t : S(\F)-\medg_{X} \longrightarrow \D_\F^{\mathrm{fil}}-\fdg_X,$$
obtained from $|-|^t$ and observe that, by definition, $\D_\F^{\mathrm{fil}}:=|S(\F)|^t$. By 
Proposition \ref{ptate}, this $\s$-functor is fully faithful and its
image consists of complete filtered $\D_\F^{\mathrm{fil}}$-modules over $X$. We restrict this
to $S(\F)-\medg_{qcoh,w=0}$, the full sub-$\s$-category of graded mixed
module which are quasi-coherent and of weight $0$. Its image by $|-|^t$ is easily seen to consist of 
all filtered $\D_\F^{\mathrm{fil}}$-modules $E$ satisfying the following two conditions

\begin{enumerate}
\item The filtration on $E$ is tautological: $F^i(E)=E$ if $i\geq 0$ and $0$ if $i<0$.
\item The underlying $\OO_X$-module of $E$ is quasi-coherent.
\end{enumerate}

Now, these two conditions define a full sub-$\s$-category of $\D_\F^{\mathrm{fil}}-\fdg_X$
which is equivalent, via the underlying object 
$\s$-functor $(-)^u : \fdg_X \to \dg_X$, to the $\s$-category
$\D_\F-\dg^{qcoh}_X$, of \emph{unfiltered} quasi-coherent $\D_\F$-modules. \hfill $\Box$ \\

The following notation will be used throughout the rest of the paper.

\begin{df}\label{dqcoh}
\emph{Let $\F \in \Fol(X)$ be a derived foliation over a derived scheme $X$. We denote by 
$\QCoh^{\mathrm{fil}}(\F)$ the full sub-$\s$-category of $\D_\F^{\mathrm{fil}}-\fdg_X$ consisting of
all filtered modules which are quasi-coherent as filtered $\OO_X$-modules via restriction of scalars along 
the natural
morphism of filtered dg-algebras $\OO_X \longrightarrow \D_\F^{\mathrm{fil}}$. We will call 
$\QCoh^{\mathrm{fil}}(\F)$  
the \emph{$\s$-category of filtered quasi-coherent crystals along $\F$}.}
\end{df}

\begin{rmk}\label{defnormcone} \emph{An interesting consequence of Theorem  
\ref{pdmod} is the existence of a \emph{deformation to the normal cone} for any derived foliation $\F \in \Fol(X)$. Indeed, 
let $\QCoh^{\mathrm{fil}}(\F)$ be the $\s$-category of quasi-coherent filtered $\D_\F^{\mathrm{fil}}$-modules over $X$.
Being the $\s$-category of filtered modules over a filtered dg-algebra, this $\s$-category 
possesses a natural tensored and cotensored structure over $\fdg\simeq\QCoh(\A=[\mathbb{A}^1/\Gm])$ (see \S \ref{sec-geompic}), the symmetric
monoidal $\s$-category of filtered complexes. Recall (\S \ref{sec-geompic}) the two 
symmetric monoidal $\s$-functors
$$(-)^u : \QCoh(\A) \longrightarrow \dg \qquad
\mathrm{Gr}(-) : \QCoh(\A) \longrightarrow \grdg.$$
It is easy to see that the underlying object and associated graded $\s$-functors on 
$\D_\F^{\mathrm{fil}}-\dg_X^{\mathrm{fil},\mathrm{qcoh}}$
induces natural equivalences
$$\QCoh^{\mathrm{fil}}(\F)\otimes_{\QCoh(\A)}\dg \simeq \QCoh(\F)$$
$$\QCoh^{\mathrm{fil}}(\F)\otimes_{\QCoh(\A)}\grdg \simeq \QCoh^{\mathrm{gr}}(\F^{\epsilon=0}).$$
Here, $\F^{\epsilon=0}$ is the derived foliation defined by 
$\DR(\F^{\epsilon=0})=Sym_{\OO_X}(\LL_{\F}[1])$ with trivial mixed structure, 
thus it is endowed with a natural $\Gm$-action. The $\s$-category $\QCoh^{\mathrm{gr}}(\F^{\epsilon=0})$
consists here of $\Gm$-equivariant quasi-coherent crystals along $\F^{\epsilon=0}$, 
and is therefore equivalent to the category of graded $Sym_{\OO_X}(\mathbb{T}_\F)$-modules over $X$. \\
Therefore, $\QCoh^{\mathrm{fil}}(\F)$ provides a family of $\s$-categories over the stack 
$\A=[\mathbb{A}^1/\Gm]$, 
whose generic fiber is $\QCoh(\F)$ and with special fiber $\QCoh(\F^{\epsilon=0})$. This
family is the deformation to the normal cone of $\F$, and degenerates $\F$ into
an abelian derived foliation. This family can also be 
constructed as a \emph{relative derived foliation on $X\times \A$ over $\A$}, 
whose graded mixed cdga is $\DR(\F)[t]$, where $t$ is the canonical parameter on $\A$ and
the mixed structure is taken to be $t.\epsilon$, where $\epsilon$ is the mixed
structure of $\DR(\F)$.
This canonical family is the starting point of 
a \emph{non-abelian Hodge theory} \`a la Simpson for derived foliations, and related notions, such as 
Higgs structures, and $\lambda$-connections along the
leaves etc. We hope to come back to this interesting subject in a later work.} 
\end{rmk}

\begin{rmk}\label{notgaro}\emph{Coming back to $X=\mathrm{Spec}\, k[u]$, where $\deg \,u =-2$, as in Remark \ref{weyl2}, Theorem \ref{pdmod} tells us that $\QCoh(*_X)$ is very different from the category $\mathrm{Crys}(X)$ of \cite{garo}. The interested reader will find in \cite{be} more informations in this direction. }
\end{rmk}

\subsection{The induction $\s$-functor}\label{s-ind}

Let $u : \F \to \G$ be a morphism of derived foliations on a derived scheme $X$. We have seen that it 
induces a morphism of filtered dg-algebras over $X$
$$\D_\F^{\mathrm{fil}} \longrightarrow \D_\G^{\mathrm{fil}}.$$
Associated to this, we have  the usual forgetful and base-change adjunction
\begin{equation}\label{for-filt}
u_!:=\D_\G^{\mathrm{fil}}\otimes_{\D_\F^{\mathrm{fil}}}- : \D_\F^{\mathrm{fil}}-\fdg_X \rightleftarrows \D_\G^{\mathrm{fil}}-\fdg_X : u^!.\end{equation}
The right adjoint $u^!$ is called the \emph{inverse image $\s$-functor}. 
The left adjoint $u_!$ is called the \emph{induction along $u$} or \emph{direct image $\s$-functor}. 
By forgetting the filtrations, we have a corresponding non-filtered adjunction
\begin{equation}\label{for-nofilt} u_!:=\D_\G\otimes_{\D_\F}- : \D_\F-\dg_X \rightleftarrows \D_\G-\dg_X : u^!.
\end{equation}
Both $\s$-functors $u_!$ and $u^!$ preserve quasi-coherence, and thus induce
an adjunction on quasi-coherent modules. By Theorem \ref{pdmod}, this can also be interpreted
as an adjunction on the $\s$-category of quasi-coherent crystals
$$u_! : \QCoh(\F) \rightleftarrows \QCoh(\G) : u^!$$
where again, $u^!$ is called the inverse image $\s$-functor, and $u_!$ the induction or direct image
$\s$-functor. By Definition \ref{dqcoh}, the filtered adjunction (\ref{for-filt}) can be regarded as an adjunction on filtered crystals, as well
$$u_! : \QCoh^{\mathrm{fil}}(\F) \rightleftarrows \QCoh^{\mathrm{fil}}(\G) : u^!.$$
The filtered and unfiltered versions of $u_{!}$ and $u^{!}$ are of course compatible with the underlying object $\s$-functor, i.e. the
following squares naturally commutes
\begin{equation}\label{dunno}\xymatrix{
\QCoh^{\mathrm{fil}}(\F) \ar[r]^{u_!} \ar[d]_-{(-)^u} & \QCoh^{\mathrm{fil}}(\G) \ar[d]^-{(-)^u} &  & \QCoh^{\mathrm{fil}}(\G) 
\ar[r]^{u^!} \ar[d]_-{(-)^u} & \QCoh^{\mathrm{fil}}(\F) \ar[d]^-{(-)^u} \\
\QCoh(\F) \ar[r]_-{u_!} & \QCoh(\G) &  & \QCoh(\G) \ar[r]_-{u^!} & \QCoh(\F)
}
\end{equation}
The same is true when the underlying object $\s$-functor $(-)^u$ is replaced with the associated graded $\s$-functor $\mathrm{Gr}$. 

\begin{rmk}
\emph{Tracking back the equivalence of Theorem \ref{pdmod} it is easy to see that the inverse image functor $u^! : \QCoh(\F) 
\longrightarrow \QCoh(\G)$
may also be identified  with the base change at the level of graded mixed dg-modules
$$\DR(\F) \otimes_{\DR(\G)} - : \DR(\G)-\medg_X \longrightarrow \DR(\F)-\medg_X$$
for the morphism $\DR(\G) \longrightarrow \DR(\F)$ corresponding to $u : \F \to \G$ in $\Fol(X)$.
In other words, it does coincide with the pull-back of quasi-coherent crystals defined in 
\S \ref{sectionqcoh} along the morphism of pairs $(id_X,u) : (X,\F) \to (X,\G)$.}
\end{rmk}

We will now examine two specific important cases of adjunctions (\ref{for-filt}) and (\ref{for-nofilt}): when either $\F$ is the initial foliation, or $\G$ is the final foliation (so that the morphism $u$ is uniquely defined in either cases).\\

Let us first consider the morphism $u: 0_X \to \F$. We know that 
$\QCoh(0_X)$ is naturally equivalent to $\QCoh(X)$. The corresponding induction $\s$-functor
$$u_! : \QCoh(X) \longrightarrow \QCoh(\F)$$
is then simply induced by $\D_\F \otimes_{\OO_X}-$. We warn the reader that 
this $\s$-functor does not have an easy description on the level of graded mixed modules, and this shows a particular instance of the usefulness of Theorem \ref{pdmod}. For example, $u_!$ 
sends $\OO_X$ to $\D_\F$ which is a rather big and complicated object inside $\QCoh(\F)$, not 
concentrated in degree $0$ (except if $\F$ and $X$ are both assumed to be smooth). The $\s$-functor
$u_!$ will play an important role for us later and will be referred to as \emph{the induction $\s$-functor
for $\F$}. There is also a corresponding filtered version 
$$u_! : \QCoh^{\mathrm{fil}}(X) \longrightarrow \QCoh^{\mathrm{fil}}(\F).$$

\begin{df}\label{dinduc}
Let $\F \in \Fol(X)$ be a derived foliation and $u : 0_X \to \F$ the canonical morphism. 
The \emph{induction $\s$-functor for $\F$} is the $\s$-functor
$$\mathsf{Ind}_\F:=u_! : \QCoh(X) \longrightarrow \QCoh(\F).$$
The \emph{filtered induction $\s$-functor for $\F$} is the $\s$-functor
$$\mathsf{Ind}^{\mathrm{fil}}_\F:=u_! : \QCoh^{\mathrm{fil}}(X) \longrightarrow \QCoh^{\mathrm{fil}}(\F).$$
\end{df}

A direct consequence of the existence of the induction $\s$-functor is the following important
observation.

\begin{cor}\label{ccompact}
The $\s$-categories $\QCoh^{\mathrm{fil}}(\F)$ and $\QCoh(\F)$ have compact generators.
\end{cor}

\noindent \textit{Proof.} We know that $\QCoh(X)$ has a  compact generator, because all our
derived schemes are assumed to be quasi-compact and quasi-separated. Pick a compact generator
$E \in \QCoh(X)$. It is formal to check that $\mathsf{Ind}_\F(E) \in \QCoh(\F)$ is a compact generator.
In the filtered case, the argument is similar, noticing that $\QCoh^{\mathrm{fil}}(X)\simeq \QCoh(X\times \A)$. 
Indeed, as $X\times \A$ is a global quotient stack of $X \times \mathbb{A}^1$ by $\Gm$ it is
again true that $\QCoh(X\times \A)$ possesses compact generators. \hfill $\Box$ \\

The other interesting special case of (\ref{for-nofilt}) is that of the unique morphism $u : \F \to *_X$ to the final foliation $*_X$. We know that $\D_{*_X}\simeq \D_X$ is the sheaf of differential operators on $X$. The induction $\s$-functor, in this
situation, produces a $\s$-functor
$$u_ ! : \QCoh(\F) \longrightarrow \D_X-\dg_X^{qcoh},$$
from quasi-coherent crystals along $\F$ to quasi-coherent $\D_X$-modules on $X$. The $\D_X$-modules 
of the form $u_!(E)$ will be called \emph{induced from the foliation $\F$}. One important example  
is the induced $\D_X$-module $u_!(\OO_X)$. This is a canonical $\D_X$-module on $X$ 
associated to the derived foliation $\F$ which contains interesting informations about $\F$. For instance, 
when $\F$ is smooth, then $u_!(\OO_X)$ is a coherent $\D_X$-module. However, this is not true anymore
for non-smooth derived foliations $\F$. Measuring the defect of coherence of $u_!(\D_X)$ is
a very interesting question related to invariants of singularities of derived foliations, generalizing 
classical invariants such as Milnor numbers.
More details about this construction will appear in a future work.

\section{Inverse and direct images of filtered crystals}
In this Section we define (filtered and unfiltered) direct image functors between quasi-coherent crystals, along proper and quasi-smooth maps. We also prove the important result that filtered direct and filtered inverse images commutes with the underlying or the associated graded objects functors.

\subsection{Direct images}\label{s-dirim} Let $f : (X,\F) \longrightarrow (Y,\G)$ be a morphism of
derived schemes endowed with derived foliations. Thus $f$ is given by a morphism
$g : X \rightarrow Y$, and a morphism of foliations on $X$
$\alpha : \F \rightarrow g^*(\G)$ (i.e. a morphism $g^*\DR(\G) \to \DR(\F)$
of graded mixed dg-algebras over $X$, see \cite{tvRH}).

We have seen in \S \ref{sectionqcoh} that $(X,\F) \mapsto \QCoh(\F)$ is a contraviant
$\s$-functor by using pull-backs. This functoriality can be extended to the
filtered case as follows. Let $f=(g,u) : (X,\F) \longrightarrow (Y,\G)$
be a morphism of pairs consisting of derived schemes and derived foliations. 
In order to define a filtered pull-back
$$f^! : \QCoh^{\mathrm{fil}}(\G) \longrightarrow \QCoh^{\mathrm{fil}}(\F),$$
we consider 
$$\D_{\G \to \F}^{\mathrm{fil}}:=g^*(\D_\G^{\mathrm{fil}})=\OO_X \otimes_{g^{-1}(\OO_Y)}g^{-1}(\D_\G^{\mathrm{fil}}).$$ 
This is a filtered $(\D_\F^{\mathrm{fil}},g^{-1}(\D_\G^{\mathrm{fil}}))$-bi-module over $X$. As such it defines a pull-back
$\s$-functor on filtered modules
$$f^! : \D_{\G}^{\mathrm{fil}}-\fdg_Y \longrightarrow \D_\F^{\mathrm{fil}}-\fdg_X$$
by  
$$f^!(E):=\D_{\G\to\F}^{\mathrm{fil}}\otimes_{g^{-1}(\D_\G^{\mathrm{fil}})}g^{-1}(E).$$
On underlying $\OO$-modules, the $\s$-functor $f^!$ acts as the usual pull-back
of filtered $\OO$-modules. In particular, 
it preserves the $\OO$-quasi-coherence conditon, and thus induces
a well defined $\s$-functor
$$f^! : \QCoh^{\mathrm{fil}}(\G) \longrightarrow \QCoh^{\mathrm{fil}}(\F).$$

\begin{lem}\label{ldirect}
With the notations above, if $g : X \to Y$ is a proper and quasi-smooth, then 
the $\s$-functor
$$f^! : \QCoh^{\mathrm{fil}}(\G) \longrightarrow \QCoh^{\mathrm{fil}}(\F)$$
admits a left adjoint 
$$f_! : \QCoh^{\mathrm{fil}}(\F) \longrightarrow \QCoh^{\mathrm{fil}}(\G).$$
\end{lem}

\textit{Proof.} We already know from Corollary \ref{ccompact} that 
both $\s$-categories are compactly generated. For the existence of 
$f_!$ it is thus enough to check that 
$f^!$ commutes with limits. For this, we use that $f^!$ is compatible
with the usual pull-backs of filtered $\OO$-modules along $g$: the following square
naturally commutes
$$\xymatrix{
\QCoh^{\mathrm{fil}}(\G) \ar[r] \ar[d]_-{f^!} & \QCoh^{\mathrm{fil}}(Y) \ar[d]^-{(g^*)^{\mathrm{fil}}} \\
\QCoh^{\mathrm{fil}}(\F) \ar[r] & \QCoh^{\mathrm{fil}}(X).}$$
Here the horizontal $\s$-functors are the forgetful functors, induced
from the natural maps of filtered dg-algebras $\OO_X \to \D_\F^{\mathrm{fil}}$
and $\OO_Y \to\D_\G^{\mathrm{fil}}$, while the $\s$-functor $(g^*)^{\mathrm{fil}}$ is the  
pull-back of filtered quasi-coherent complexes, obtained by 
applying $\mathrm{Fun}(\ZZ^{\leq},-)$ to the usual pull-back $g^* : \QCoh(Y) \longrightarrow \QCoh(X)$. 
The horizontal $\s$-functors are clearly conservative and commute with 
limits and colimits. Therefore, to check that $f^!$ 
preserves limits, it is enough to show that $(g^*)^{\mathrm{fil}}$ does. Again, 
as $(g^*)^{\mathrm{fil}}$ is obtained from the usual pull-back $g^* : \QCoh(Y) \longrightarrow \QCoh(X)$
by applying $\mathrm{Fun}(\ZZ^{\leq},-)$, we are reduced to show that $g^*$ preserves limits. 

This last step follows easily from the assumption that $g$ is proper and quasi-smooth. Indeed, 
let $E \in \QCoh(X)$ be a compact generator (thus a perfect 
complex on $X$), and $\{F_i\}_{i\in I}$ a diagram in $\QCoh(Y)$. 
Using that $E$ is dualizable and projection formula, we have
$$Map(E,g^*(\lim_i F_i)) \simeq Map(\OO_X,E^\vee \otimes_{\OO_X} g^*(\lim_i F_i)) \simeq
Map(\OO_Y,g_*(E^\vee)\otimes_{\OO_Y} (\lim_i F_i)).$$
Now we use that $g$ is proper and quasi-smooth, so that $g_*(E^\vee)$ is again perfect and thus dualizable, 
and therefore the functor $g_*(E^\vee)\otimes_{\OO_Y} -$ commutes with limits. So we have
$$Map(\OO_Y,g_*(E^\vee)\otimes_{\OO_Y} (\lim_i F_i)) \simeq \lim_i Map(\OO_Y,g_*(E^\vee)\otimes_{\OO_Y} F_i) \simeq
\lim_i Map(E,g^*(F_i)).$$
This shows that the canonical
map
$$Map(E,g^*(\lim_i F_i)) \longrightarrow \lim_i Map(E,g^*(F_i))\simeq Map(E,\lim_i \, g^*(F_i))$$
is indeed in equivalence, so that $g^*$ preserves limits (since $E$ is a generator of $\QCoh(Y)$). \hfill $\Box$ \\

Lemma \ref{ldirect} allows us to give the following

\begin{df}\label{ddirect}
Let $f=(g,u) : (X,\F) \longrightarrow (Y,\G)$ be a morphism with 
$g : X \to Y$ proper and quasi-smooth. The \emph{filtered direct image} is the $\s$-functor 
$$f_! : \QCoh^{\mathrm{fil}}(\F) \longrightarrow \QCoh^{\mathrm{fil}}(\G),$$
right adjoint to the pull-back $\s$-functor $f^!$.
\end{df}

\begin{rmk}
\emph{Note that, although it is possible, we do not try to define direct images 
for non-proper or non-quasi-smooth morphisms. In the rest of the paper, we will only need this kind of direct images.}
\end{rmk}

We also define the unfiltered direct image
$$f_! : \QCoh(\F) \longrightarrow \QCoh(\G)$$
as being the left adjoint to the unfiltered version of $f^!$. \\

The direct image functors satisfy the usual pseudo-functoriality properties, $(ff')_!\simeq f_!f_!'$ (for composeable, proper and quasi-smooth $f$ and $f'$). An important 
consequence of this property is the following result.

\begin{prop}\label{cdirectinduct}
Let $f=(g,u) : (X,\F) \longrightarrow (Y,\G)$ be a morphism with $g$ proper and quasi-smooth.
Then, 
the following diagram naturally commutes
$$\xymatrix{
\QCoh^{\mathrm{fil}}(X) \ar[r]^-{\mathsf{Ind}^{\mathrm{fil}}_\F} \ar[d]_-{g_*(\omega_{X/Y}\otimes -)[d]} & 
\QCoh^{\mathrm{fil}}(\F) \ar[d]^-{f_!} \\
\QCoh^{\mathrm{fil}}(Y) \ar[r]_-{\mathsf{Ind}^{\mathrm{fil}}_\G} & \QCoh^{\mathrm{fil}}(\G),
}$$
where $\omega_{X/Y}$ is the relative canonical line bundle of $X$ over $Y$ and 
$d$ the relative dimension of $X$ over $Y$.
\end{prop}

\textit{Proof.} We have a commutative diagram of pairs
$$\xymatrix{
(X,0_X) \ar[r]^-{u} \ar[d]_-{g} & (X,\F) \ar[d]^-{f} \\
(Y,0_Y) \ar[r]_-{v} & (Y,\G)
}$$
where $u$ and $v$ are the unique morphism from the initial foliation. We get from this a natural isomorphism
of $\s$-functors 
$$f_!u_!\simeq v_!g_!.$$
The proposition then follows from the fact that, by definition $u_!$ and $v_!$ are the induction $\s$-functors, 
and from the explicit formula for $g_! : \QCoh(X) \longrightarrow \QCoh(Y)$
in terms of relative Serre duality. \hfill $\Box$ \\

An interesting application of Proposition \ref{cdirectinduct} is to the direct image of $\D_\F^{\mathrm{fil}}$ itself. With the same notations as in Proposition \ref{cdirectinduct}, we obtain 
$$f_!(\D_\F^{\mathrm{fil}}) \simeq \D_\G^{\mathrm{fil}}\otimes_{\OO_Y}g_*(\omega_{X/Y})[d].$$

When $X$ is proper and quasi-smooth over $k$,
and $f : (X,\F) \longrightarrow (Spec\, k,0)$ is the projection to the point (endowed with its trivial
foliation), 
we see in particular that $f_!(\D_\F)$ computes $H^{*+d}(X,\omega_X)$, that is \emph{coherent homology 
of $X$ with coefficients in $\OO_X$.} \\

The following result is a direct consequence of the fact that 
$f^!$ commutes with colimits (and can also be deduced  
by  using Proposition \ref{cdirectinduct}
and the proof of Corollary \ref{ccompact}). 

\begin{cor}\label{cproper}
Let $f=(g,u) : (X,\F) \to (Y,\G)$ be a morphism of derived schemes endowed with derived
foliations, with $g$ proper and quasi-smooth. Then, the direct image $\s$-functors
$$f_! : \QCoh^{\mathrm{fil}}(\F) \longrightarrow \QCoh^{\mathrm{fil}}(\G) \qquad 
f_! : \QCoh(\F) \longrightarrow \QCoh(\G)$$
preserve compact objects.
\end{cor}

\subsection{Compatibility with underlying and associated graded objects} We conclude this Section with the important result that 
filtered direct images commute with both taking underlying and associated graded objects.
For this, let us consider a morphism $f=(g,u) : (X,\F) \longrightarrow (Y,\G)$ 
with $g$ proper and quasi-smooth. We have the corresponding  adjunction on filtered crystals
$$f^{\mathrm{fil}}_! : \QCoh^{\mathrm{fil}}(\F) \rightleftarrows \QCoh^{\mathrm{fil}}(\G) : f_{\mathrm{fil}}^! \, ,$$
and its unfiltered version
$$f_! : \QCoh(\F) \rightleftarrows \QCoh(\G) : f^!.$$
We may also consider $\F^{\epsilon=0}$ which is the derived foliation 
whose underlying graded mixed cdga is $\DR(\F)^{\epsilon=0}$ (i.e. with trivial mixed structure $\epsilon=0$), 
and the same for $\G^{\epsilon=0}$. For such derived foliations, we have now a graded push-forward
$$f^{\mathrm{gr}}_! : \QCoh^{\mathrm{gr}}(\F^{\epsilon=0}) \longrightarrow \QCoh^{\mathrm{gr}}(\G^{\epsilon=0})$$
defined as the left adjoint to the graded pull-back $f^!_{\mathrm{gr}} : \QCoh^{\mathrm{gr}}(\G^{\epsilon=0}) \rightarrow 
\QCoh^{\mathrm{gr}}(\F^{\epsilon=0})$.
Note that the $\s$-categories $\QCoh^{\mathrm{gr}}(\F^{\epsilon=0})$ and $\QCoh^{\mathrm{gr}}(\G^{\epsilon=0})$
can also be identified with the categories of graded $Sym_{\OO_X}(\T_\F)$-modules over $X$ and of 
graded $Sym_{\OO_Y}(\T_\G)$-modules over $Y$, respectively (where $\T_\F$ and $\T_\G$ both sit in weight $1$).

By putting all these functors together, we may write the following diagram of vertical adjunctions
\begin{equation}\label{compGr}
\xymatrix{
\QCoh^{\mathrm{gr}}(\F^{\epsilon=0}) \ar@<-1.0ex>[d]_-{f^{\mathrm{gr}}_!} & \QCoh^{\mathrm{fil}}(\F) \ar[l]_-{\mathrm{Gr}} \ar[r]^-{(-)^u}
\ar@<-1.0ex>[d]_-{f^{\mathrm{fil}}_!} & \QCoh(\F) \ar@<-1.0ex>[d]_-{f_!} \\
\QCoh^{\mathrm{gr}}(\G^{\epsilon=0}) \ar@<-1.0ex>[u]_-{f_{\mathrm{gr}}^!} & 
\QCoh^{\mathrm{fil}}(\G) \ar[l]^-{\mathrm{Gr}} \ar[r]_-{(-)^u}
\ar@<-1.0ex>[u]_-{f_{\mathrm{fil}}^!} & \QCoh(\G) \ar@<-1.0ex>[u]_-{f^!}
}
\end{equation}

We already know that this diagram, when restricted to the inverse images only, naturally commutes. This
implies that the diagram restricted to direct images is naturally lax commutative. In fact, the natural
transformations
$$f_!^{\mathrm{gr}}\circ \mathrm{Gr} \Rightarrow \mathrm{Gr}\circ f^{\mathrm{fil}}_! 
\qquad
f_!\circ (-)^u \Rightarrow (-)^u\circ f_!^{\mathrm{fil}}$$
turn out to be equivalences. Indeed, as all $\s$-functors involved commute with
colimits, it is enough to check this property on compact generators, and  we can thus use
Proposition \ref{cdirectinduct}
and the proof of Corollary \ref{ccompact} to conclude. As a consequence of the commutativity of (\ref{compGr}), we are allowed to (and will from now on) simply
write $f_!$ and $f^!$ \emph{without} any decorations $(-)^{\mathrm{fil}}$ or $(-)^{\mathrm{gr}}$. 
Because of its importance, and for later reference, we state this result in the following corollary.

\begin{cor}\label{ccommute}
For a morphism $f=(g,u) : (X,\F) \longrightarrow (Y,\G)$ with $g$ proper and quasi-smooth, 
filtered direct and filtered inverse images of quasi-coherent crystals commute with 
taking the underlying or the associated graded objects, i.e. the following diagram
naturally commutes.
\begin{equation}
\xymatrix{
\QCoh^{\mathrm{gr}}(\F^{\epsilon=0}) \ar@<-1.0ex>[d]_-{f_!} & \QCoh^{\mathrm{fil}}(\F) \ar[l]_-{\mathrm{Gr}} 
\ar[r]^-{(-)^u}
\ar@<-1.0ex>[d]_-{f_!} & \QCoh(\F) \ar@<-1.0ex>[d]_-{f_!} \\
\QCoh^{\mathrm{gr}}(\G^{\epsilon=0}) \ar@<-1.0ex>[u]_-{f^!} & 
\QCoh^{\mathrm{fil}}(\G) \ar[l]^-{\mathrm{Gr}} \ar[r]_-{(-)^u}
\ar@<-1.0ex>[u]_-{f^!} & \QCoh(\G) \ar@<-1.0ex>[u]_-{f^!}
}
\end{equation}
\end{cor}

\begin{rmk}\emph{Here is an equivalent way of looking at the leftmost adjunction $f^{\mathrm{gr}}_! :\QCoh^{\mathrm{gr}}(\F^{\epsilon=0}) \leftrightarrows \QCoh^{\mathrm{gr}}(\G^{\epsilon=0}) : f^{\mathrm{gr}\, !}$ of (\ref{compGr}) that avoids introducing the auxiliary foliations $\F^{\epsilon=0}$ and $\G^{\epsilon=0}$. First of all, we have the associated graded object $\s$-functor 
$$\mathrm{Gr}: \QCoh^{\mathrm{fil}}(\F)= \mathrm{Mod}_{\D_\F^{\mathrm{fil}}}(\fdg_X)^{\mathrm{qcoh \, on \,  }X} \longrightarrow \mathrm{Mod}_{\mathrm{Gr}(\D_\F^{\mathrm{fil}})}(\grdg_X)^{\mathrm{qcoh \, on \,  }X},$$ and the analogous one for $\G$ and $Y$. Note that  $\mathrm{Gr}(\D_\F^{\mathrm{fil}}) \simeq Sym_{\OO_X}(\mathbb{T}_\F)$ and $\mathrm{Gr}(\D_\G^{\mathrm{fil}}) \simeq Sym_{\OO_Y}(\mathbb{T}_\G)$ in 
$\mathrm{Alg}(\grdg_X)$ (with $\mathbb{T}_\F$ and $\mathbb{T}_\G$ of pure weight $1$). 
Now, we proceed as in \S \, \ref{s-dirim}. Define
$$\mathrm{Gr}(\D^{\mathrm{fil}})_{\G \to \F}:=g^*(\mathrm{Gr}(\D_\G^{\mathrm{fil}}))=\OO_X \otimes_{g^{-1}(\OO_Y)}g^{-1}(\mathrm{Gr}(\D_\G^{\mathrm{fil}})).$$ 
which is a graded $(\mathrm{Gr}(\D_\F^{\mathrm{fil}}),g^{-1}(\D_\G^{\mathrm{fil}}))$-bi-module over $X$, and as such, it defines a pull-back
$\s$-functor on graded modules
$$f_{\mathrm{gr}}^! : \mathrm{Mod}_{\mathrm{Gr}(\D_{\G}^{\mathrm{fil}}))}(\grdg_Y) \longrightarrow \mathrm{Mod}_{\mathrm{Gr}(\D_\F^{\mathrm{fil}})}(\grdg_X)$$
by  
$$f_{\mathrm{gr}}^!(E):=\mathrm{Gr}(\D^{\mathrm{fil}})_{\G \to \F}\otimes_{g^{-1}(\mathrm{Gr}(\D_\G^{\mathrm{fil}}))}g^{-1}(E).$$
This $\s$-functor respects the property of being quasi-coherent over $X$ and $Y$, so it induces a $\s$-functor
$$f_{\mathrm{gr}}^! : \mathrm{Mod}_{\mathrm{Gr}(\D_{\G}^{\mathrm{fil}}))}(\grdg_Y)^{\mathrm{qcoh \, on \,  }Y} \longrightarrow \mathrm{Mod}_{\mathrm{Gr}(\D_\F^{\mathrm{fil}})}(\grdg_X)^{\mathrm{qcoh \, on \,  }X}$$ which has a left adjoint 
$$f^{\mathrm{gr}}_!:\mathrm{Mod}_{\mathrm{Gr}(\D_\F^{\mathrm{fil}})}(\grdg_X)^{\mathrm{qcoh \, on \,  }X} \leftrightarrows \mathrm{Mod}_{\mathrm{Gr}(\D_\G^{\mathrm{fil}})}(\grdg_Y)^{\mathrm{qcoh \, on \,  }X} : f_{\mathrm{gr}}^!,$$ which coincides with the leftmost adjunction  of (\ref{compGr}), since $\QCoh^{\mathrm{gr}}(\F^{\epsilon=0}) \simeq \mathrm{Mod}_{\mathrm{Gr}(\D_\F^{\mathrm{fil}})}(\grdg_X)^{\mathrm{qcoh \, on \,  }X}$, and $\QCoh^{\mathrm{gr}}(\G^{\epsilon=0}) \simeq \mathrm{Mod}_{\mathrm{Gr}(\D_\G^{\mathrm{fil}})}(\grdg_Y)^{\mathrm{qcoh \, on \,  }Y}$.
}
\end{rmk}

\section{Characteristic cycles}

In this section we introduce the notion of \emph{characteristic cycle} of
a quasi-coherent crystal (or $\D_\F$-modules thanks to the equivalence
of Theorem \ref{pdmod}) along a derived foliation $\F$. For this we first introduce the global cotangent stack $T^*\F$ of a derived
foliation $\F$, which is a derived Artin n-stack, where $n$ is the tor-amplitude of the perfect complex
$\LL_\F$. We will then discuss
the notion of \emph{bounded
coherent crystals} and of \emph{good filtrations} on them. By definition, the associated graded 
to a good filtration will be a $\Gm$-equivariant perfect complex on $T^*\F$, that will be used to define 
 characteristic cycles. We investigate the 
existence of good filtrations and prove some independence (of good filtrations) results for 
characteristic cycles.

\subsection{Cotangent stacks of derived foliations}

Let $X$ be a derived scheme and $E$ be a perfect complex on $X$ of amplitude 
contained in $[0,n]$ for some non-negative integer $n$ (see
\cite{moduli} for the notion of amplitude of perfect complexes). We can associate to $E$ 
a linear stack $\VV(E)$ over $X$, whose functor of points sends $u : Spec\, A \to X$ to 
the space $Map(u^*(E),A)$, of morphisms of $A$-modules from $u^*(E)$ to $A$. As shown in 
\cite[Sub-lemma 3.9]{moduli} $\VV(E)$ is a smooth Artin $n$-stack over $X$. Moreover, the
projection $\pi : \VV(E) \longrightarrow X$ makes it into a linear derived stack over $X$.
The stack $\VV(E)$ comes equipped with an obvious 
$\Gm$-action, covering the morphism $\pi$, by acting on $E$ via its natural weight $1$ action. 
This makes makes $\VV(E)$ into a $\Gm$-equivariant derived Artin stack over $X$. \\

We consider
the morphism on quotient stacks
$$\pi^{\mathrm{gr}} : [\VV(E)/\Gm] \longrightarrow X \times B\Gm.$$
The direct image along this morphism is a symmetric lax monoidal $\s$-functor
\begin{equation}\label{for-this}\pi_*^{\mathrm{gr}} : \QCoh^{\mathrm{gr}}(\VV(E)):=\QCoh([\VV(E)/\Gm]) \longrightarrow \QCoh(X \times B\Gm)=\QCoh^{\mathrm{gr}}(X).\end{equation}
Since structure sheaf $\OO:= \OO_{[\VV(E)/\Gm]}$ is the monoidal unit in $\QCoh^{\mathrm{gr}}(\VV(E))$, the lax-monoidal the $\s$-functor $\pi_*^{\mathrm{gr}}$ in (\ref{for-this}) factors via a $\s$-functor (denoted by the same symbol)
$$\pi_*^{\mathrm{gr}} : \QCoh^{\mathrm{gr}}(\VV(E)) \longrightarrow \pi^{\mathrm{gr}}_*(\OO)-\grdg_{qcoh,X},$$
from graded quasi-coherent complexes on $\VV(E)$ to 
graded $\pi^{\mathrm{gr}}_*(\OO)$-modules over $X$ which are quasi-coherent as $\OO_X$-modules.

\begin{prop}\label{plinear}
There exists a fully faithful $\s$-functor
$$Sym_{\OO_X}(E)-\grdg_{\mathrm{perf},X} \longrightarrow \QCoh^{\mathrm{gr}}(\VV(E))$$
from 
perfect graded $Sym_{\OO_X}(E)$-modules over $X$, sending
the i-th twist $Sym_{\OO_X}(E)(i)$ of the tautological graded module $Sym_{\OO_X}(E)$, to the i-th twist  $\OO_{\VV(E)}(i)$ of the structure sheaf of $\VV(E)$.
\end{prop}

\textit{Proof.} We start by assuming that $X=Spec\, A$ is affine, with $A$ a
connective cdga, and we construct  
an equivalence of graded cdga's over $X$
$$\pi^{\mathrm{gr}}_*(\OO)\simeq Sym_{\OO_X}(E).$$
We consider $E^\vee$, the dual perfect $A$-dg-module of $E$, which is perfect of amplitude
$[-n,0]$. By the Dold-Kan equivalence, we can write $E$ as the colimit of a simplicial 
diagram $n \mapsto E'_n$ of vector bundles over $Spec\, A$. As the functor
of $p$-th symmetric power commutes with sifted colimits, we have, for all $p\geq 0$, a natural equivalence
$$\mathrm{colim}_{n\in \Delta^o}Sym^{p}_A(E'_n) \simeq Sym_A^p(E^\vee)$$ of perfect $A$-modules.
This is a colimit inside the $\s$-category of perfect $A$-modules, so we can dualize these equivalence
to get equivalences
$$Sym^p_A(E^\vee)^{\vee} \simeq Sym^p_A(E) \simeq \lim_{n\in \Delta}Sym_A^p(E_n),$$
where $E_n$ is the dual $A$-module of $E_n'$. Taking the sum over all $p$, we get an equivalence of graded
$A$-modules
$$Sym_A(E) \simeq \lim_{n \in \Delta}Sym_A(E_n).$$
By construction, this is clearly an equivalence of graded $A$-linear cdga's. \\

We now consider the simplicial diagram of derived stacks $n \mapsto \VV(E_n)$. By construction,
the natural morphism $\mathrm{colim}_n\VV(E_n) \longrightarrow \VV(E)$ is an equivalence of $\Gm$-equivariant
derived stacks over $X$. We thus find an induced equivalence of graded $A$-linear cdga's
$$\OO^{\mathrm{gr}}(\VV(E)) \simeq \lim_n \OO^{\mathrm{gr}}(\VV(E_n)),$$
where $\OO^{\mathrm{gr}}$ denotes the graded cdga of functions on $\Gm$-equivariant derived stacks
(taking values in possibly non-connective graded cdga's). But 
each $E_n$ is a vector bundle, and thus $\VV(E_n)=Spec\, Sym_A(E_n)$ is an affine derived stack
over $X$. We thus have a natural equivalence of $A$-linear graded cdga's
$Sym_A(E_n) \simeq \OO^{\mathrm{gr}}(\VV(E_n))$, functorial in $n \in \Delta$. Assembling these facts together, 
we get the required equivalence of $A$-linear graded cdga's
$$u : \OO^{\mathrm{gr}}(\VV(E)) \simeq Sym_A(E).$$
The equivalence $u$ is clearly functorial in $A$, and thus can be globalized over a more general
derived scheme $X$. We thus find a natural equivalence of quasi-coherent graded $\OO_X$-linear cdga's
over $X$
$$\pi_*^{\mathrm{gr}}(\OO) \simeq Sym_{\OO_X}(E).$$

We are now ready to conclude the proof of the proposition. We have already considered the natural $\s$-functor
$$\pi_*^{\mathrm{gr}} : \QCoh^{\mathrm{gr}}(\VV(E)) \longrightarrow \pi^{\mathrm{gr}}_*(\OO)-\grdg_{qcoh,X},$$
from graded quasi-coherent complexes on $\VV(E)$ to 
graded $\pi_*(\OO)$-modules over $X$ which are quasi-coherent as $\OO_X$-modules. When 
restricted to the thick triangulated sub-$\s$-category generated by the objects
$\OO(i)$, this induces an equivalences with perfect graded $\pi^{\mathrm{gr}}_*(\OO)$-modules over $X$. 
As $\pi^{\mathrm{gr}}_*(\OO)\simeq Sym_{\OO_X}(E)$, 
the inverse of this equivalence is the $\s$-functor in the proposition.
\hfill $\Box$ \\

The previous proposition will be applied later to the particular
case where $E=\T_\F$, the tangent complex of a derived foliation $\F$ on $X$. 
The derived stack $V(\T_\F)$ will be denoted by $T^*\F$, and will be called
the \emph{global derived cotangent stack} of $\F$.

\begin{df}\label{dcotangent}
For a derived scheme $X$ and a derived foliation $\F \in \Fol(X)$, the \emph{(derived) cotangent stack 
of $\F$} is defined by 
$$T^*\F:=V(\T_\F).$$
equipped with is natural $\Gm$-action.
\end{df}

Proposition \ref{plinear} will be used to define characteristic cycles, 
by applying the $\s$-functor in the proposition to the associated graded of filtered quasi-coherent crystals along $\F$.

\subsection{Good filtrations on crystals}

Let $X$ be a derived scheme and $\F \in \Fol(X)$ be a derived foliation on $X$. We let 
$E \in \QCoh(\F)$ be a quasi-coherent crystal. Via the equivalence of Theorem 
\ref{pdmod}, we will freely identify 
$E$ with a $\D_\F$-module quasi-coherent over $X$.

\begin{df}\label{dgood}
Let $X$, $\F$ and $E$ as above. 
\begin{enumerate}
\item The crystal $E$ is called \emph{coherent} if it is a compact object
in $\QCoh(\F)$. The full sub-$\s$-category of coherent crystals is denoted by 
$\Coh(\F) \subset \QCoh(\F)$.

\item For $E \in \Coh(\F)$, a \emph{good filtration on $E$} is the data of
a compact object $E^{\mathrm{fil}} \in \QCoh^{\mathrm{fil}}(\F)$ together with an equivalence in $\QCoh(\F)$
$$(E^{\mathrm{fil}})^u \simeq E.$$
\end{enumerate}
\end{df}

Compact objects in $\QCoh(\F)$ and in $\QCoh^{\mathrm{fil}}(\F)$ can be easily characterized, either
via the induction $\s$-functors, or as sheaves of $\D_\F$-modules.

\begin{prop}\label{pcompact}
Let $X$ and $\F$ be as above. An object $E \in \QCoh(\F)$ (resp.  
$E^{\mathrm{fil}} \in \QCoh^{\mathrm{fil}}(\F)$) is compact if and only if it
satisfies one of the following two equivalent conditions.

\begin{enumerate}
\item The object $E$ belongs to the thick triangulated sub-category generated
by objects of the form $\mathsf{Ind}_{\F}(E_0)$ (resp. $\mathsf{Ind}^{\mathrm{fil}}_{\F}(E^{\mathrm{fil}}_0)$)
for a perfect complex of $\OO_X$-modules $E_0$ (resp. a filtered perfect
complex of $\OO_X$-modules $E_0^{\mathrm{fil}}$).

\item $E$ (resp. $E^{\mathrm{fil}}$) is a perfect $\D_\F$-module (resp. a perfect $\D_\F^{\mathrm{fil}}$-module) i.e. locally on $X_{zar}$, $E$ (resp. $E^{\mathrm{fil}}$) is
a retract of a finite cell $\D_\F$-module (resp. of  a finite cell filtered $\D_\F^{\mathrm{fil}}$-module).

\end{enumerate}
\end{prop}

\textit{Proof.} Condition $(1)$ have been already considered in the proof of Corollary \ref{ccompact}. 
Clearly, condition $(1)$ implies condition $(2)$. Finally condition $(2)$ clearly implies 
compactness when $X$ is affine. The general case follows from
the quasi-compactness of $X$. \hfill $\Box$ \\

Suppose that $E$ is a coherent crystal along $\F$, equipped with a good filtration 
$E^{\mathrm{fil}}$ in the sense of the definition above. By Proposition \ref{pcompact}, 
$\mathrm{Gr}(E^{\mathrm{fil}})$ is a perfect graded $\mathrm{Gr}(\D_\F^{\mathrm{fil}})=Sym_{\OO_X}(\T_\F)$-module over $X$. 
Now, Proposition \ref{plinear} implies that $\mathrm{Gr}(E^{\mathrm{fil}})$ defines a graded perfect 
complex on the stack $T^*\F$, and by forgetting the $\Gm$-action 
we get a perfect complex on $T^*\F$. We consider the $K$-theory spectrum
$K(T^*\F)$, defined as the $K$-theory of the $\s$-category of perfect
(not graded) modules over $T^*\F$. The perfect complex $\mathrm{Gr}(E^{\mathrm{fil}})$
on $T^*\F$ thus defines a class
$$Ch(E^{\mathrm{fil}}) \in K_0(T^*\F).$$

\begin{df}\label{dsing1}
The \emph{characteristic cycle} of $E^{\mathrm{fil}}$ is the element
$$Ch(E^{\mathrm{fil}}) \in K_0(T^*\F)$$
defined above.
\end{df}

Note that the above definition depends a priori on the choice of $E^{\mathrm{fil}}$. We will see later
that, in fact, modulo phantoms, it does not (see \S \, \ref{sect-inde}, and Proposition \ref{pinv}). \\

\subsubsection{Existence of good filtrations.} The existence of good filtrations in general
seems a complicated question, and the authors do not know if good filtrations always
exist for coherent crystals, as it is the case for usual $\D$-modules on smooth varieties. 
It can be shown that they do exist for smooth foliations on smooth varieties, but the fact that, for general $\F$, its ring of differential operators $\D_\F$ is not concentrated
in degree $0$, creates complications in constructing good filtrations. The following result is therefore very 
useful in practice.

\begin{prop}\label{pgooddirect}
Let $f=(g,u) : (X,\F) \longrightarrow (Y,\G)$ be a morphism of smooth varieties endowed with
derived foliations and assume that $g$ is proper. If $E^{\mathrm{fil}}$ is a good filtration
on a coherent crystal $E \in \QCoh(\F)$, then $f_!(E^{\mathrm{fil}})$ is a good
filtration on $f_!(E)$.
\end{prop}

\textit{Proof.} This follows easily from Corollary \ref{cproper} and from the fact that 
direct images commutes with taking the underlying object (see diagram (\ref{dunno})). \hfill $\Box$ \\

We can also isolate a large class of coherent crystals for which existence
of good filtrations is guaranteed: the finite cell crystals. For
a derived scheme $X$ and $\F \in \Fol(X)$, we define a non-thick 
triangulated sub-$\s$-category $\Coh^{cell}(\F) \subset \Coh(\F)$, as being generated
(by finite limits and shifts)
by the objects of the form $Ind_\F(E)$ for $E$ a perfect complex on $X$. An object 
of $\Coh^{cell}(\F)$ will be called  a \emph{finite cell crystal}. There is an obvious filtered
version, too.
Finite cell crystals are
obviously coherent, but we do not know if the converse is true. Finite cell crystals are
however useful because of the following result.

\begin{prop}\label{pcell}
With the notation above, any object $E \in \Coh^{cell}(\F)$ admits 
a good filtration. Moreover, a good filtration $E^{\mathrm{fil}}$ can be chosen to
be a finite cell object in $\QCoh^{\mathrm{fil}}(\F)$.
\end{prop}

\textit{Proof.} By definition of being a finite cell object, there is a finite sequence of morphisms
in $\Coh(\F)$
$$\xymatrix{
E_{-1}=0 \ar[r] & E_0 \ar[r] & \dots \ar[r] & E_i \ar[r] & E_{i+1} \ar[r] & \dots \ar[r] & 
E_n=E,}$$
with the following property: for all $i$ there is a perfect complex $K_i$ on 
$X$ and a cartesian square in $\Coh(\F)$
$$\xymatrix{
E_i \ar[r] & E_{i+1} \\
Ind_\F(K_i) \ar[u]^-{u_i} \ar[r] & 0.  \ar[u]}$$
We can show, by induction, that $E_i$ has a good filtration. For this, assume
that $E_i$ has a good filtration $E_i^{\mathrm{fil}}$. The morphism $u_i$ is given, by adjunction, by a morphism 
$v_i : K_i \longrightarrow E_i=(E_i^{\mathrm{fil}})^u$ in $\QCoh(\F)$.
The quasi-coherent sheaf $(E_i^{\mathrm{fil}})^u$ is the filtered colimit $\mathrm{colim}_k F^k(E_i^{\mathrm{fil}})$, and
as $K_i$ is a compact object in $\QCoh(X)$, $v_i$ can be factored as
$\xymatrix{
K_i \ar[r]^-{w_i} & F^k(E_i^{\mathrm{fil}}) \ar[r] & E_i
},$
for some index $k$. By using the left adjoint
$Ind_\F^{\mathrm{fil}}$, the morphism $w_i$ corresponds to a morphism of filtered $\D_\F^{\mathrm{fil}}$-modules
$$\alpha_i :  Ind_{\F}^{\mathrm{fil}}(K_i)(-k) \longrightarrow E_i,$$
where $(-k)$ denotes the endofunctor of $\QCoh^{\mathrm{fil}}(\F)$ that shifts by $-k$ the filtration. The cone 
of $\alpha_i$ clearly defines a good filtration on $E_{i+1}$.
\hfill $\Box$ \\

\subsubsection{Independence of the good filtration.}\label{sect-inde} We already noticed that the characteristic 
cycle of Definition 
\ref{dsing1} \emph{depends}, a priori, on the good filtration $E^{\mathrm{fil}}$. In order to solve this
problem, we introduce a reduced $K$-group, and prove that the image of $Ch(E^{\mathrm{fil}})$ in this reduced $K$-group, only
depends on the object $E$. \\

Let $\F$ be a derived foliation on a derived scheme $X$. 
We consider $\Coh^{\mathrm{fil}}(\F)$, the $\s$-category of compact objects $\QCoh^{\mathrm{fil}}(\F)$, and
the underlying object $\s$-functor $(-)^u : \Coh^{\mathrm{fil}}(\F) \longrightarrow \Coh(\F)$.
An object $E \in \Coh^{\mathrm{fil}}(\F)$ will be called a \emph{phantom} if 
$E^u\simeq 0$. We then set the following definition.

\begin{df}\label{reduced}
With the notations above, the reduced $K$-group $K^{red}_0(T^*\F)$ is the 
quotient of $K_0(T^*\F)$ by the subgroup generated by the classes of $Gr(E)$ for 
$E \in \Coh^{\mathrm{fil}}(\F)$ a phantom.
\end{df}

We then have the following result.

\begin{prop}\label{pinv}
Let $\F$ be a derived foliation on a derived scheme $X$, and
$E \in \QCoh(\F)$ be a coherent crystal.
Let $E_1^{\mathrm{fil}}$ and $E_2^{\mathrm{fil}}$
be two good filtrations on $E$ (in the sense of Definition \ref{dgood}) which are finite cell
filtered object.
Then we have $Ch(E_1^{\mathrm{fil}})=Ch(E_2^{\mathrm{fil}})$ in $K^{red}_0(T^*\F)$.
\end{prop}

\textit{Proof.} We start by the following lifting lemma.

\begin{lem}\label{llift}
Let $M,N \in \QCoh^{\mathrm{fil}}(\F)$ with $M$ compact. Then 
any morphism $u : M^u \longrightarrow 
N^u$ in $\QCoh(\F)$ can be lifted, via the $\s$-functor $(-)^u$, to a morphism
$v : M \longrightarrow N(k)$ to the $k$-shift of $N$, for some integer $k$.
\end{lem}

\textit{Proof of the lemma.} Recall that $N(k)$ denotes the filtered crystal 
with the filtration shifted by $k$, i.e.
$$F^i(N(k)):=F^{i+k}N.$$
As $M$ is compact, we know by the proof of Corollary \ref{ccompact}, that 
$M$ is a retract of a finite cell object in $\QCoh^{\mathrm{fil}}(\F)$ in the sense of \ref{pcell}. 
Clearly, if the lemma is true for $M$ (and any $N$) it is also true for any of its retracts. We 
may therefore assume that $M$ is a finite cell object. By induction on the number of cells we 
reduce to the following statement. Assume that the lemma is true for $M$ (and any $N$), and let 
us consider a push-out
$$\xymatrix{
M \ar[r] & M' \\
Ind_{\F}(K) \ar[u]-^{\alpha} \ar[r] & 0 \ar[u]
}$$
with $K$ compact in $\QCoh(X)$.
We must prove that the lemma remains true for $M'$ (and any $N$). Let $v : (M')^u \longrightarrow
N^u$ be a morphism. It consists of the data of a morphism $u : M^u \longrightarrow N^u$
and a homotopy to zero of $u\alpha : Ind_{\F}(K) \longrightarrow N^u$, or by adjunction 
a homotopy to zero $h$ of the induced morphism $\beta : K \longrightarrow N^u$. 
By assumption on $u$ we can lift $u$ to $w : M \longrightarrow N(k)$ for some $k$. Moreover, 
 $h$ defines a homotopy to zero of $v=(w)^u : M^u \longrightarrow N(k)^u \simeq N^u$. 
As $K$ is compact and $N^u=\mathrm{colim}_k F^iN$, the homotopy $h$ factors as a homotopy to zero
$h'$ of $K \longrightarrow F^{k'}N \longrightarrow N^u$ for some $k'\geq k$. This pair $(w,h')$
defines the required lift $M' \longrightarrow N(k')$. 
\hfill $\Box$ \\

Let us go back to the proof of Proposition \ref{pinv}. As $E_i^{\mathrm{fil}}$, $i=1, 2$ are a good filtrations on the same crystal $E$, we have a natural 
equivalence in $\QCoh(F)$
$$u : (E_1^{\mathrm{fil}})^u \simeq (E_2^{\mathrm{fil}})^u.$$
By Lemma \ref{llift}, $u$ can be lifted to a morphism of filtered crystals
$v : E_1^{\mathrm{fil}} \longrightarrow E_2^{\mathrm{fil}}$. We set 
$M$ to be the cone of $v$ in $\QCoh^{\mathrm{fil}}(\F)$. This is a compact object which is obviously a
phantom, i.e. $M^u\simeq 0$. We thus have an $Ch(E_1^{\mathrm{fil}})=Ch(E_2^{\mathrm{fil}})$ in $K^{red}_0(T^*\F)$.
\hfill $\Box$ \\

Proposition \ref{pinv} shows that the following notion is well defined.

\begin{df}
Let $X$ be a derived scheme and $\F \in \Fol(X)$. If $E\in \Coh(\F)$ admits a good filtration $E^{\mathrm{fil}}$, 
then its characteristic cycle is defined as
$$Ch(E)=Ch(E^{\mathrm{fil}}) \in K^{red}_0(T^*\F).$$
\end{df}

It is possible to show that when $X$ and $\F$ are both smooth, then 
the natural projection $K_0(T^*\F) \longrightarrow K_0^{red}(T^*\F)$
is bijective, or, in other words, that the class of $Gr(E)$ is trivial
in $K_0(T^*\F)$ for any phantom $E \in \Coh^{\mathrm{fil}}(\F)$. This relies on 
using regularity and Quillen's devissage techniques, that do not work in 
our general setting. We will not need this isomorphism $K_0(T^*\F) \simeq K_0^{red}(T^*\F)$ in the rest of the paper, so we omit its proof . 
However, the following very simple particular case will be useful later,
in order to get numerical formulas out of our general Grothendieck-Riemann-Roch statement.

\begin{prop}\label{ppoint}
Let $X$ be a derived scheme and $\F=0_X$ be the initial foliation
so that $T^*\F\simeq X$. Then, the natural projection
$$K_0(X) \longrightarrow K_0^{red}(X)$$
is bijective.
\end{prop}

\textit{Proof.} Let $\QCoh^{\mathrm{fil}}_0(X) \subset \QCoh^{\mathrm{fil}}(X)$ be the full sub-$\s$-category 
defined to be the kernel of the underlying object $\s$-functor $(-)^u$. 
By using the relation between graded mixed complexes and filtered objects 
(Proposition \ref{ptate}), we see that $\QCoh^{\mathrm{fil}}_0(X)$ is equivalent $\medg_{X,qcoh}$, the
$\s$-category of graded mixed quasi-coherent complexes on $X$. Forgetting the graded structures
yields an $\s$-functor
$$\QCoh^{\mathrm{fil}}_0(X) \simeq \medg_{X,qcoh} \longrightarrow \edg_{X,qcoh},$$
to non-graded mixed quasi-coherent complexes on $X$. By definition of mixed complexes, the
$\s$-category on the right hand side can be natural identified with 
$\QCoh(X[\epsilon])$, where $X(\epsilon):= =X\times Spec\, k[\epsilon]$ with $\deg \, \epsilon= -1$. Passing to compact objects, we get 
$$\phi : \Coh^{\mathrm{fil}}_0(X) \longrightarrow \Parf(X[\epsilon]).$$
Now, the map
$Gr : K_0(\Coh^{\mathrm{fil}}_0(X)) \longrightarrow K_0(X)$
clearly factors as
$$Gr : \xymatrix{
K_0(\Coh^{\mathrm{fil}}_0(X)) \ar[r]^-{\phi} & K_0(X[\epsilon]) \ar[r]^{p_*} & K_0(X),
}$$
where $p : X[\epsilon]=X\times Spec\, k[\epsilon] \longrightarrow X$ is the first projection.
We are thus reduced to show that $p_*$ is zero in $K$-theory. The projection $p$
has a section $j : X \to X[\epsilon]$, and thus $j^*p^*=id$, so that  $p^* : K_0(X) \to K_0(X[\epsilon])$ is an injective map. It is thus enough to prove that 
$p^*p_*$ is the zero endomorphism of $K_0(X[\epsilon])$. But, clearly, for a perfect complex
$E$ on $X[\epsilon]$, we have 
$$p^*p_*(E)\simeq E\oplus E[1].$$
This shows that $p^*p_*$ is zero in $K$-theory, and implies the proposition.
\hfill $\Box$ \\

A particularly important case of Proposition \ref{ppoint} is when $X=\mathrm{Spec}\, k$, for which
we find $K_0^{red}(\mathrm{Spec}\, k)\simeq \ZZ$.

\section{GRR for derived foliations}

In this Section we will state and prove the Grothendieck-Riemann-Roch formula for proper maps between
derived schemes endowed with derived foliations. 

\subsection{The GRR formula}

Let 
$f=(g,u) : (X,\F) \longrightarrow (Y,\G)$ be a morphism of derived schemes with 
derived foliaions, with $g$ proper and quasi-smooth.
Associated to $f$ is the 
so-called ``Japanese correspondence''
$$\xymatrix{
T^*\F & T^*\G \times_Y X \ar[r]^-{p} \ar[l]_-{q} & T^*\G.}$$
Here $q$ is induced by the morphism of $g^*(\LL_\G) \to \LL_\F$ perfect complexes on $X$ 
induced by $u$. The morphism $p$ simply is the first projection.

Define the push-forward on K-groups
$$f_!:=p_!(q^*(-)) : 
K_0(T^*\F) \longrightarrow K_0(T^*\G).$$
Here $p_!$ is the left adjoint of the pull-back of quasi-coherent sheaves
$p^*  : \QCoh(T^*\G\times_Y X) \longrightarrow \QCoh(T^*\G)$. This left adjoint exists 
because $p$ is representable, proper and quasi-smooth, and it is given explicitly by 
$$p_!(E)=p_*(E\otimes \omega_p[d])$$ for $E\in \Parf(T^*\G\times_Y X)$.
Here, the integer
$d$ is the relative dimension of $X$ over $Y$, and $\omega_p$ is the relative canonical sheaf of the 
morphism $p$, which is also the
pull-back of $\omega_{X/Y}$ along the projection $T^*\G\times_Y X \to X$. 

We start by noticing that $f_!$ is compatible with the quotient defining 
reduced $K$-groups. Indeed, if $E \in \Coh^{\mathrm{fil}}(\F)$ is a phantom, then so is $f_!(E)$, and
the image by $f_!$ of $Gr(E)$ is $Gr(f_!E)$
(see Corollary \ref{ccommute}). Therefor, $f_{!}$ induces a well defined map
$$f_! : K_0^{red}(T^*\F) \longrightarrow K_0^{red}(T^*\G).$$

\begin{thm}\label{tgrr}
Let $f=(g,u): (X,\F) \longrightarrow (Y,\G)$ be a morphism of 
derived schemes endowed with derived foliations. If $g$ is proper and quasi-smooth, and 
$E \in \Coh(\F)$ is a coherent crystal along $\F$ admitting a good filtration (e.g. 
a finite cell object), then we have
$$Ch(f_!(E))=f_!(Ch(E))$$
as elements of $K_0^{red}(T^*\G)$.
\end{thm}

\textit{Proof.} This is a direct consequence of the fact that direct images
commutes with taking associated graded for good filtrations (Corollary \ref{ccommute}). \hfill $\Box$ \\

The following corollary of Theorem \ref{tgrr} is obtained when $Y=\mathrm{Spec}\, k$ and $\G=0_Y$. 
Let $s : X \to T^*\F$ denotes the zero section of the
cotangent stack. When $X$ is a proper and quasi-smooth derived scheme we denote by 
$p : X \to Spec\, k$ the canonical map, and by
$$p_*=:\int_X : K_0(X) \longrightarrow K_0(k)\simeq \ZZ$$
the push-forward on perfect complexes.

\begin{cor}\label{chrr}
Let $(X,\F)$ be a quasi-smooth and proper derived scheme endowed with a derived foliation $\F\in \Fol(X)$. 
Let $f : (X,\F) \to (Spec\, k,0_{\mathrm{Spec}\, k})$ be the canonical morphism.
For any coherent crystal $E$ that admits a good filtration we have
$$\chi(f_!(E))=\int_{X}s^*(Ch(E))\otimes \omega_X.$$
\end{cor}

This corollary is a \emph{Hirzeburch-Riemann-Roch} (HRR) \emph{formula} for crystals along the foliation $\F$. The complex
$f_!(E)$ is what should be called the \emph{foliated cohomology of $(X,\F)$ with coefficients in $E$}.
If we denote this cohomology by $H^*_{\F}(X,E)$, the HRR formula reads
$$\chi(H^*_{\F}(X,E))=\int_{X}s^*(Ch(E))\otimes \omega_X.$$

\subsection{The non-proper case: Fredholm crystals}\label{s-ell}

We explain here briefly how to extend Theorem \ref{tgrr} to non-proper maps, 
by introducing the notion of \emph{Fredholm crystal}. \\

We start by extending direct images to the case of \emph{compactifiable morphisms}. 
Assume that $f=(g,u) : (X,\F) \to (Y,\G)$ is a morphism of derived schemes endowed with derived foliations, and $g: X \to Y$ is quasi-smooth. A \emph{quasi-smooth compactification} of $f$ is the datum of a factorization
\begin{equation}\label{fact}
f : \xymatrix{
(X,\F) \ar[r]^-{(j,v)} & (\bar{X},\bar{\F}) \ar[r]^-{\bar{f}} & (Y,\G)}
\end{equation}
where $j$ is an open embedding, $v$ is an equivalence $\F \simeq j^*(\bar{\F})$, 
and $p$ is a proper and quasi-smooth smooth morphism. For such a compactification, 
we set 
$$j_* : \QCoh(\F) \to \QCoh(\bar{\F})$$
to be the right adjoint to the $\s$-functor $j^!$. We note here that $j_*$ exists because
$j^!$ commutes with colimits. Moreover, $j_*$ is compatible with the usual push-forward of quasi-coherent
sheaves by the forgeful $\s$-functors, i.e. the following diagram naturally commutes
$$\xymatrix{
\QCoh(\F) \ar[r]^-{j_*} \ar[d] & \QCoh(\bar{\F}) \ar[d] \\
\QCoh(X) \ar[r]_-{j_*} & \QCoh(\bar{X}),
}$$
where the vertical $\s$-functors are the natural forgetful functors. Because 
$j$ is an open immersion and $v$ is an equivalence, we clearly have a natural 
equivalence of sheaves of dg-algebras on $X$
$j^{-1}\D_{\bar{\F}} \simeq \D_\F$, so that we get a canonical adjunction 
morphism on $\bar{X}$ 
$$a : \D_{\bar{\F}} \longrightarrow j_*(\D_\F).$$
The $\s$-functor $j_*$ is then simply induced by the usual push-forward of 
sheaves along $j : X \hookrightarrow \bar{X}$: it sends a
$\D_\F$-module $E$ to $j_*(E)$, viewed as a $\D_{\bar{\F}}$-module via the map
$a$ above. This description shows that we also have a commutative square 
involving the induction $\s$-functor
$$\xymatrix{
\QCoh(\F) \ar[r]^-{j_*}  & \QCoh(\bar{\F}) \ \\
\QCoh(X) \ar[r]_-{j_*} \ar[u]^-{Ind_{\F}} & \QCoh(\bar{X}). \ar[u]_-{Ind_{\bar{\F}}}
}$$

We define the \emph{$*$-direct image along $f$} as
$$f_*=\bar{f}_!\circ j_* : \QCoh(\F) \longrightarrow \QCoh(\G).$$
A standard argument, by using the product embedding, proves that $f_*$, defined as above, does not depend on the
choice of the factorization (\ref{fact}) (see e.g. \cite[Exp XVII]{SGA4.3} or \cite[\S 8]{FrKi}). 
Indeed, $f_*$ clearly commutes with colimits, so 
the independence of the factorization in its definition can be checked on compact $\D_\F$-modules
of the form $Ind_\F(E)$, where it is easily deduced from the compatibility of push-forwards
for quasi-coherent sheaves, and from the explicit formula of Proposition \ref{cdirectinduct}. 

\begin{df}\label{dell}
Let $f=(g,u) : (X,\F) \to (Y,\G)$ be a morphism of derived schemes endowed with derived foliations. Suppose that $f$ admits a quasi-smooth compactification. An object $E \in \Coh(\F)$ is called \emph{f-Fredholm}
if it admits a good filtration $E^{\mathrm{fil}} \in \Coh^{\mathrm{fil}}(\F)$ such that $j_*(E^{\mathrm{fil}})$ is a compact object in $\QCoh^{\mathrm{fil}}(\F)$ for some quasi-smooth compactification $j : (X,\F) \hookrightarrow
(\bar{X},\bar{\F})$ of the morphism $f$.
\end{df}

Suppose $E \in \Coh(\F)$ is $f$-Fredholm; pick a quasi-smooth compactification $(\bar{X},\bar{\F})$ 
and a filtration $E^{\mathrm{fil}}$ as in the definition above. 
We have the associated graded $Gr(E^{\mathrm{fil}}) \in \Parf(T^*\F)$, and its direct image by $j$ is by 
again a perfect complex on $T^*\bar{\F}$, so that, in particular, the support 
of $Gr(E^{\mathrm{fil}})$, $SS(E) \subset T^*\F$ is closed in $T^*\bar{\F}$. If we pull-back 
this support by the canonical map $T^*\G \times_Y X \to \T^*\F$,  
we get a closed subset in $T^*\G \times_Y X$ which is proper over $T^*\G$. 
In other words, the morphism $f$ is proper when restricted to the support of $Gr(E^{\mathrm{fil}})$. 
%This condition is the classical \emph{ellipticity condition} for $\D$-modules (see for instance ????). 
In our situation, the notion of Fredholm is a priory stronger, it is unclear to us that 
properness of $f$ on the support of $Gr(E^{\mathrm{fil}})$ is enough to recover the fact that 
$j_*(E^{\mathrm{fil}})$ remains a compact object. Therefore, our definition of being $f$-Fredholm might
be hard to check in practice. \\

We can now state the Grothendieck-Riemann-Roch formula for possibly non-proper maps, 
and Fredholm coefficients. 

\begin{thm}\label{tgrrell}
Let $f=(g,u): (X,\F) \longrightarrow (Y,\G)$ be a morphism of 
derived schemes endowed with derived foliations. If $f$ admits a quasi-smooth
compactification, and $E \in \Coh(\F)$ is an $f$-Fredholm crystal,
then we have an equality 
$$Ch(f_*(E))=f_*(Ch(E))$$
of elements in $K_0^{red}(T^*\G)$.
\end{thm}

\textit{Proof.} Simply apply Theorem \ref{tgrr} to the morphism $\bar{f}$ and to the 
object $j_*(E)$ over $(\bar{X},\bar{\F})$,
which has a good filtration given by $j_*(E^{\mathrm{fil}})$ for $E^{\mathrm{fil}}$ given by 
Definition \ref{dell}. \hfill $\Box$ \\

\section{Examples and applications}

We finish this paper by giving a sample of examples and applications of the general Grothendieck-Riemann-Roch formula of Theorem \ref{tgrr}. \\

\subsection{Grothendieck-Riemann-Roch for derived $\D$-modules.}\label{s-grrdmod}

The very first special case of the theorem \ref{tgrr} is when $\F=*_X$ and $\G=*_Y$ are the
final derived foliations on $X$ and $Y$ respectively. The $\s$-categories
$\QCoh(*_X)$ and $\QCoh(*_Y)$ are then called the \emph{$\s$-categories of derived 
$\D$-modules on $X$ and $Y$}, and denoted by $\D_X-\dg_{X,qcoh}$ and $\D_Y-\dg_{Y,qcoh}$, respectively. 
We think they coincide with the derived $\D$-modules introduced and studied in \cite{be}, and denoted there by $\D^{der}(X)$ and $\D^{der}(Y)$. While it is clear that $\D^{der}(X) \simeq \D_X-\dg_{X,qcoh}$ for $X= \mathrm{Spec} (k[u])$, with $\deg (u)= -n$, for arbitrary $n\geq 0$, we will not attempt to give here a precise general comparison between our categories and the ones considered in \cite{be}. However, we remark that, contrary to the notion of crystals introduced in \cite{garo}, these $\s$-categories of derived $\D$-modules (both the one in \cite{be} and ours) are sensitive to derived structures (see also Remark \ref{notgaro}).

In this setting, $T^*\F=T^*X$ and $T^*\G=T^*Y$ are the global derived cotangent stacks of $X$ and $Y$, and the GRR formula is the equality  
\begin{equation}\label{ddmodform}Ch(f_!(E))=f_!Ch(E) \,\,\, \textrm{in} \,\, K_0^{red}(T^*Y) \end{equation}
for a compact object $E \in \D_X-\dg_{X,qcoh}$ admitting a good filtration, and $f : X \to Y$ 
a proper and quasi-smooth morphism.\\ When $X$ and $Y$ are smooth varieties, this recovers\footnote{Note that, as already observed, for any smooth variety $S$ over $k$, the foliation $*_S=\DR (S)$ is smooth, and the canonical map $K_0(T^*S) \to K_0^{red}(T^*S)$ is bijective. So our result exactly recovers the usual GRR for $\D$-modules on smooth varieties. } the well known 
GRR formula for $\D$-modules of \cite{la} (see also \cite{sab}). However, already when $X$ and $Y$ are
just underived $k$-schemes, formula (\ref{ddmodform}) is new and provides a Grothendieck-Riemann-Roch formula for $\D$-modules on \emph{possibly singular} schemes. \\

An interesting feature of this situation is that any morphism $f=(g,u) : (X,\F) \to (Y,\G)$ enters with in a commutative
diagram
$$\xymatrix{
(X,\F) \ar[r]^-{f} \ar[d]_-{u} & (Y,\G) \ar[d]^-{v} \\
(X,*_X) \ar[r]_-g & (Y,*_Y).
}$$
As explained at the end of \S \, \ref{s-ind}, the push-fowards $u_!$ and $v_!$ define
induction $\s$-functors \emph{from} the given foliations 
$$u_! : \QCoh(\F) \to \D_X-\dg_{X,qcoh} \qquad
v_! : \QCoh(\G) \to \D_Y-\dg_{Y,qcoh}.$$
The GRR formula for $f$, which is an equality in $K^{red}(T^*\G)$, can be
pushed foward along $v$ to a formula in $K_0^{red}(T^*Y)$:
$$v_!Ch(f_!(E))=g_!Ch(u_!E))$$
where $E$ is a coherent crystal along $\F$ admitting a good filtration, 
$u_!(E)$ the induced coherent $\D_X$-module, and $v_! : K^{red}_0(T^*\G) \to K^{red}_0(T^*Y)$
the direct image along the canonical morphism of foliations $\G \to *_Y$. This gives a formula for the push-forward of
the characteristic cycle of a $\D_X$-module \emph{induced from} the derived foliation $\F$. When $f=id$ the formula reads
$$v_!Ch(E)=Ch(u_!E)$$
and should be understood as a formula for the characteristic cycle on a $\D_X$-module induced
from the derived foliation $\F$. \\

\subsection{Smooth Lie algebroids} Let $X$ be a smooth variety and $\F \in \Fol(X)$ a smooth foliation on $X$ (i.e. $\LL_{\F}$ is a vector bundle on $X$). As explained in 
\cite{tvRH}, $\F$ is defined by a Lie algebroid $\T_\F \to \T_X$. The sheaf $\D^{\mathrm{fil}}_\F$ is then isomorphic to the universal
enveloping algebra $U(\T_\F)$ of $\T_\F$, equipped with its natural PBW filtration. 
As already noticed $K_0^{red}(T^*\F)\simeq K_0(T^*\F)$, in this case (of a smooth foliation on a smooth variety). It can also be shown that any coherent crystal along $\F$
admits a good filtration.
Moreover, $\D_\F$ satisfies all the 
conditions of Quillen theorem (\cite{quillen}) and we thus have natural isomorphisms of K-groups
$$\tau_X : K_0(\Coh(\F)) \simeq K_0(T^*\F) \simeq K_0(X).$$
The first of this isomorphisms is precisely given by $E \mapsto Gr(E^{\mathrm{fil}})$ for $E^{\mathrm{fil}}$ a good filtration on $E$, and
the second isomorphism is the pull-back along the zero section $s : X \to T^*\F$.

The GRR formula \ref{tgrr} tells us that the isomorphism $\tau_X : K_0(\Coh(F)) \simeq K_0(X)$ is 
covariantly functorial in $(X,\F)$
in the following sense: given a morphism $f=(g,u) : (X,\F) \to (Y,\G)$ of smooth varieties endowed
with smooth foliations, with $g$ proper, the following diagram
$$\xymatrix{
K_0(\Coh(\F)) \ar[r]^-{\tau_X} \ar[d]_-{f_!} & K_0(X) \ar[d]^-{g_*(-\otimes \omega_{X/Y})[d]} \\
K_0(\Coh(\G)) \ar[r]_-{\tau_Y} & K_0(Y)
}$$
commutes. This recovers the well known GRR formula for $\D$-modules on smooth varieties, and its natural extension
to Lie algebroids. This extension to Lie algebroids is probably a folklore result (as its proof is word by word the same as for $\D$-modules), but we have not been able to find a reference
in the literature. \\

\subsection{Shifted Poisson structures} An important class of derived foliations are given by \emph{shifted Poisson structures}
in the sense of \cite{cptvv}. Indeed, let $X$ be a derived scheme endowed with a 
shifted Poisson structure of degree $n$. The Poisson bracket defines 
a morphism of perfect complexes on $X$
$$a : \LL_X[-n] \to \T_X$$
making $\LL_X[-n]$ into a dg-Lie algebroid over $X$. This defines a derived foliation $\F$ 
whose underlying graded mixed cdga is $Sym_{\OO_X}(\T_X[n])$, where the mixed structure is
induced the bracket $[-,p]$, $p$ being the bivector defining the Poisson structure. We do not know any reference
where this construction has been carried out in details, but the it can be carried out 
by representing the Poisson structure by an actual strict $\mathbb{P}_{1-n}$-structure on the
structure sheaf $\OO_X$ thanks to the strictification result of \cite{mel}.

This derived foliation associated to a shifted Poisson structure is the derived analogue of the foliation by symplectic leaves of a classical
Poisson strcuture on a smooth variety. Its \emph{leaves}, in the sense of \cite{tvRH}, are by definition
the symplectic leaves of the shifted Poisson structure. When the Poisson structure is non-degenerate, then 
the foliation $\F$ is the final foliation. In general $\F$ is a very interesting derived foliation containing 
information about the Poisson structure. For instance, our notion $\QCoh(\F)$ of quasi-coherent crystals provides a useful setting to study various versions of Poisson cohomology. \\

As an example, we may define \emph{derived Poisson cohomology} to be $$H_{der}^{\mathrm{Poiss}}(X):=\rh_{\QCoh(\F)}(\OO_X,\OO_X),$$ 
which coincides with derived de Rham cohomology when the Poisson structure is a symplectic structure. It is also
possible to consider the induced $\D_X$-module $u_!(\OO_X) \in \D_X-\dg_{X,qcoh}$, where
$u : \F \to *_X$ is the unique morphism to the final foliation. The $\D_X$-module
$u_!(\OO_X)$ can also be used in order to define another  version of Poisson homology by 
$$H_{\mathrm{Poiss}}(X):=p_*(u_!(\OO_X)) \in \dg$$
where $p : X \to Spec\, k$ is the structure map (assuming here that $X$ is quasi-smooth and  $p$ admits 
a quasi-smooth compactification as defined in \S \, \ref{s-ell}). This offers a criterion for finiteness of Poisson homology
by requiring $u_!(\OO_X)$ to be Fredholm as a $\D_X$-module on $X$. When $X$ is smooth and 
the Poisson structure is a classical Poisson structure of degree $0$, the $\D_X$-module
$u_!(\OO_X)$ has been considered in \cite{brtr} where it was used to define  the notion of \emph{holonomic Poisson varieties}
and to get finiteness results for Poisson cohomology. The GRR formula \ref{tgrrell}, and more generally 
the general formalism of crystals along derived foliations, provides a way to extend these notions and results  
to shifted Poisson structures. \\

\subsection{A foliated index formula} As a last application of theorem \ref{tgrrell} we propose an \emph{index formula} for
a weakly Fredholm differential operator along the leaves of a derived foliation on a quasi-smooth derived
scheme. We like to think of it as an algebraic version of the longitudinal index theorem of \cite{cosk}, 
possibly valid outside the smooth setting (i.e. for non-smooth derived schemes and non-smooth
derived foliation). 

Let $(X,\F)$, where $X$ is a derived scheme and $\F\in \Fol(X)$. We assume that 
$X$ has a quasi-smooth compactification $j : X \hookrightarrow \bar{X}$ (i.e.  the projection $X\to Spec\, k$ admits a quasi-smooth compactification as in \S \, \ref{s-ell}). Note that we only
assume the existence of $\bar{X}$, and we do not assume that $\F$ can be extended to $\bar{X}$.
Let $E$ and $E'$ be two perfect complexes on $X$. A \emph{differential operator $P$ along $\F$
from $E$ to $E'$} is by definition a morphism of quasi-coherent sheaves on $X$
$$P : E \longrightarrow Ind_\F(E')=\D_\F \otimes_{\OO_X}E',$$
or, equivalently, a morphism of coherent crystals along $\F$
$$a_P : Ind_\F(E) \longrightarrow Ind_\F(E').$$
As $E$ is compact in $\QCoh(X)$, we note that $P$ must factor through a morphism
$E \to \D_\F^{\leq i}\otimes_{\OO_X}E'$ for some integer $i$. In this case, we say that 
$P$ is of order $\leq i$.

We let $\D(P)$ be the cone of the morphism $a_P$ inside $\Coh(\F)$. 
%\textcolor{blue}{We say that 
%$P$ is \emph{weakly Fredholm} if the induced $\D_X$-module 
%$\D_X \otimes_{\D_\F} D(P)$ is Fredholm as a coherent crystal on $(X,*_X)$, where
%$*_X$ is the final derived foliation. \texttt{This is spurious, right?}}

We want to apply  Theorem \ref{tgrrell} to the coherent crystal $\D_X \otimes_{\D_\F} \D(P)$. Note that 
this is also the cone of the morphism $\D_X \otimes_{\D_\F} E \to \D_X \otimes_{\D_\F} E'$, induced by the composition
$$P : \xymatrix{
E \ar[r] & \D_\F \otimes_{\OO_X}E' \ar[r] & \D_X \otimes_{\OO_X}E'. }$$
In order to apply Theorem \ref{tgrrell},  we need to impose a condition on $P$ insuring that $\D_X \otimes_{\D_\F} \D(P)$ is a Fredholm object (Definition \ref{dell}).
We consider the least integer $i$ such that $P$ factors as
$$E \longrightarrow \D^{\leq i}_\F \otimes_{\OO_X}E'.$$
Associated to this, we have a morphism of filtered $\D_\F^{\mathrm{fil}}$-modules
$$Ind_\F^{\mathrm{fil}}(E)(-i) \to Ind_\F^{\mathrm{fil}}(E),$$
whose cone defines a good filtration $\D^{\mathrm{fil}}(P)$ on $\D(P)$, and thus, by 
base change, a good filtration $\D^{\mathrm{fil}}_X \otimes_{\D^{\mathrm{fil}}_\F} \D^{\mathrm{fil}}(P)$ on 
$\D_X \otimes_{\D_\F} D(P)$. 

\begin{df}\label{dellop}
The operator $P$ along $\F$ is \emph{weakly Fredholm} if 
$j_*(\D^{\mathrm{fil}}_X \otimes_{\D^{\mathrm{fil}}_\F} \D^{\mathrm{fil}}(P))$ is a compact object 
in $\QCoh^{\mathrm{fil}}(\bar{X})$, for a quasi-smooth compactification $j : X \hookrightarrow \bar{X}$.
\end{df}

Assuming that $P$ is weakly Fredholm, we can apply Theorem \ref{tgrrell}
to get our index formula, that has values in $K_0(k)\simeq \ZZ$, and thus is
an equality of two numbers we are now going to describe.\\
Let $f : X \to Spec\, k$ the projection, that we factor
as $\xymatrix{X \ar[r]^-j & \bar{X} \ar[r]^-{\bar{f}} & Spec\, k.}$
We first describe $f_*(\D(P))\simeq \bar{f_!}(j_*(\D(P)))\simeq \bar{f}_!(\D_X \otimes_{\D_\F} D(P))$.
This is a perfect complex of $k$-modules, which by definition, is the cone of the morphism induced by $P$
$$\Gamma(P) : \Gamma(X,E\otimes \omega_X)[d] \to \Gamma(X,E'\otimes \omega_X)[d].$$
Note that none of the two complexes $\Gamma(X,E\otimes \omega_X)$ or $\Gamma(X,E'\otimes \omega_X)$
is perfect, and only the cone of $\Gamma(P)$ is so. This is the effect of the weakly Fredholm property, implying that 
$\Gamma(P)$ is indeed a Fredlhom operator. Therefore, we have a first well defined number, called
the \emph{algebraic index of $P$} 
$Ind(P)$, which is the Euler characteristic of the cone of $\Gamma(P)$
$$Ind(P):=(-1)^d.\chi(cone \left(\Gamma(P) : \Gamma(X,E\otimes \omega_X) \to \Gamma(X,E'\otimes \omega_X)\right)).$$

On the other hand, the object $\D_X \otimes_{\D_\F} \D(P)$ is endowed with the good filtration
$\D^{\mathrm{fil}}_X \otimes_{\D^{\mathrm{fil}}_\F} \D^{\mathrm{fil}}(P)$, and by assumption its associated graded
defines a perfect complex of $T^*X$ which remains perfect on $T^*\bar{X}$. 
This associated graded can be described as follows. Recall that 
$i$ is the least integer such that $P$ defines a morphism in $\QCoh(X)$
$P : E \to \D_X^{\leq i} \otimes_{\OO_X}E'.$
By projection, we find 
$$\sigma(P)  : E \to Sym^i_X(\T_X) \otimes_{\OO_X}E'$$
which is called the \emph{symbol of $P$}. This extends to $$Sym_X(\T_X) \otimes_{\OO_X}E \to Sym_X(\T_X) \otimes_{\OO_X}E$$
as a morphism of graded $Sym_{\OO_X}(\T_X)$-modules 

The cone of this morphism defines a perfect complex, still denoted by $\sigma(P)$, on 
$T^*X$, which remains perfect on $T^*\bar{X}$ by the weakly Fredholm assumption.
If we denote by $s : X \to T^*X$ the zero section map, we thus get a well defined perfect complex of
$k$-modules
$\Gamma(X,s^*(\sigma(P))\otimes \omega_X)[d]$. The Euler characteristic of this complex is called
the \emph{K-theoretical index of the operator $P$}, and denoted by 
$$Ind_\sigma(P):=(-1)^d\chi(\Gamma(X,s^*(\sigma(P))\otimes \omega_X)).$$

Theorem \ref{tgrrell} implies the following

\begin{cor}\label{cindex}
With the notations above, we have
$$Ind(P) = Ind_\sigma(P).$$
\end{cor}

In plain words, Corollary \ref{cindex}, says that the \emph{algebraic index of $P$}, computing the Fredholm index of $P$ acting on global sections of $E$ and $E'$, equals 
the \emph{$K$-theoretic index of $P$}, which should be understood as an intersection number between $X$ and the cycle defined by the symbol of $P$ inside
the total cotangent stack $T^*X$.

\bibliographystyle{alpha}
\bibliography{Fol2.bib}

\end{document}